\documentclass[11pt, leqno]{amsart}
\usepackage[active]{srcltx} 
\usepackage{amsmath}
\usepackage{amssymb}
\usepackage{amsfonts}
\usepackage[usenames]{color}
\usepackage{version}

 \newtheorem{theorem}{Theorem}[section]
 \newtheorem{corollary}[theorem]{Corollary}
 \newtheorem{lemma}[theorem]{Lemma}
 \newtheorem{proposition}[theorem]{Proposition}

\newtheorem{definition}[theorem]{Definition}
\newtheorem{remark}[theorem]{Remark}
\newtheorem{example}[theorem]{Example}
\newtheorem{fact*}{Fact}

\DeclareMathOperator{\RE}{Re}
\DeclareMathOperator{\IM}{Im}

\DeclareMathOperator{\dist}{dist}

\DeclareMathOperator\re{Re}
\DeclareMathOperator\im{Im}
\newcommand\half{\tfrac 12}
\newcommand\dd{\mathrm d}

\newcommand{\s}{\mathcal{S}_n}

\newcommand\h{\mathcal{H}}

\newcommand{\M}{\mathcal{M}}
\newcommand{\N}{\mathcal{N}}
\newcommand{\T}{\mathbb{T}}

\newcommand{\D}{\mathbb{D}}
\newcommand{\C}{\mathbb{C}}

\newcommand{\pick}{\mathcal{L}_n}
\newcommand{\schur}{\mathcal{S}_n}
\newcommand{\R}{\mathbb{R}}

\newcommand{\cc}[1]{\overline{#1}}
\newcommand{\abs}[1]{\left\vert#1\right\vert}

\newcommand{\norm}[1]{\left\Vert#1\right\Vert}

\newcommand{\ran}[1]{\operatorname{ran}#1}
\newcommand{\nt}{\stackrel{\mathrm {nt}}{\to}}

\newcommand{\ip}[2]{\left\langle #1, #2 \right\rangle}
\newcommand{\ad}{^\ast}
\newcommand{\inv}{^{-1}}

\newcommand{\threepartdef}[6]
{
	\left\{
		\begin{array}{lll}
			#1 & \mbox{ if } #2 \\
			#3 & \mbox{ } #4 \\
			#5 & \mbox{ if } #6
		\end{array}
	\right.
}
\renewcommand\L{\mathcal{L}}
\newcommand\Pick{\mathcal P}
\newcommand{\p}{\mathcal{L}_n}

\newcommand{\vp}{\varphi}
\newcommand{\ph}{\varphi}
\newcommand\al{\alpha}

\newcommand\la{\lambda}

\newcommand\beq{\begin{equation}}
\newcommand\ds{\displaystyle}
\newcommand\eeq{\end{equation}}
\newcommand\df{\stackrel{\rm def}{=}}

\newcommand\black{\color{black}}
\newcommand\red{\color{red}}

\newcommand\nn{\nonumber}
\newcommand\bbm{\begin{bmatrix}}
\newcommand\ebm{\end{bmatrix}}
\numberwithin{equation}{section}
\newcommand\nin{\noindent}
\begin{document}

\title[Nevanlinna representations]{Nevanlinna representations in several variables}

\author{J. Agler, R. Tully-Doyle and N. J. Young}

\thanks{The first author was partially supported by National Science Foundation Grant on  Extending Hilbert Space Operators DMS 1068830. The third author was partially supported by the UK Engineering and Physical Sciences Research Council grant EP/J004545/1.}
\date{23rd June, 2012}
\begin{abstract}
We generalize to several variables the classical theorem of Nevanlinna that characterizes the Cauchy transforms of positive measures on the real line.  We show that for the Loewner class, a large class of analytic functions that have non-negative imaginary part on the upper polyhalfplane, there are representation formulae in terms of densely-defined self-adjoint operators on a Hilbert space.  We identify four types of such representations, and we obtain function-theoretic conditions that are necessary and sufficient  for a given function to possess a representation of each of the four types.\\
\end{abstract}
\maketitle
\section{Introduction}\label{intro}
In a classic paper \cite{nev22} of 1922 R. Nevanlinna solved the problem of the determinacy of solutions of the Stieltjes moment problem.   {\em En route} he proved several other theorems that have since been influential; in particular, the following theorem, which characterizes the Cauchy transforms of positive finite measures $\mu$ on $\R$, has had a profound impact on the development of modern analysis. Let $\Pick$ denote the Pick class, that is, the set of analytic functions on the upper halfplane,
\[
\Pi \df  \{z \in \C: \im z > 0\},
\]
that have non-negative imaginary part on $\Pi$.
 \begin{theorem}[Nevanlinna's Representation] \label{thm1.0}
  Let $h$ be a function defined on $\Pi$. There exists a finite positive measure $\mu$ on $\R$ such that
  \beq \label{cauchyt}
  h(z) = \int \frac{d\mu}{t-z}
  \eeq
  if and only if $h \in \Pick$ and
  \beq \label{cauchycond}
  \liminf_{y\to\infty} y\abs{h(iy)}< \infty.
  \eeq
 \end{theorem}

 A closely related theorem, also referred to in the literature as Nevanlinna's Representation, provides an integral representation for a general element of $\Pick$.
 \begin{theorem}\label{thm1.1}
  A function $h:\Pi\to \C$  belongs to the Pick class $ \mathcal P$ if and only if there exist $a\in\R, \ b\geq 0$ and a finite positive Borel measure $\mu$ on $\R$ such that
  \beq\label{classical}
   h(z) = a+bz+\int\frac{1+tz}{t-z} \ \dd\mu(t)
  \eeq
for all $z\in\Pi$.  Moreover, for any $h\in\Pick$, the numbers $a\in\R, \ b\geq 0$ and the measure $\mu \geq 0$ in the representation \eqref{classical} are uniquely determined.
 \end{theorem}

What are the several-variable analogs of Nevanlinna's theorems?  In this paper we shall propose four types of Nevanlinna representation for various subclasses of the $n$-variable Pick class $\Pick_n$, where $\Pick_n$ is defined to be the set of analytic functions $h$ on the polyhalfplane $\Pi^n$ such that $\im h \geq 0$. In addition, we shall present necessary and sufficient conditions for a function defined on $\Pi^n$ to possess a representation of a given type in terms of asymptotic growth conditions at $\infty$.

The integral representation \eqref{cauchyt} of those functions in the Pick class that satisfy condition \eqref{cauchycond} can be written in the form
\[
h(z) = \ip{(A-z)\inv\mathbf{1}}{\mathbf{1}}_{L^2(\mu)},
\]
where $A$ is the operation of multiplication by the independent variable on $L^2(\mu)$ and $\mathbf{1}$ is the constant function $1$.  We propose that an appropriate $n$-variable analog of the Cauchy transform is the formula
\beq\label{ctn}
h(z_1,\dots,z_n) = \ip{(A-z_1Y_1-\dots -z_nY_n)\inv v}{v}_\h \qquad \mbox{ for } z_1,\dots,z_n \in\Pi,
\eeq
where $\h$ is a Hilbert space, $A$ is a densely defined self-adjoint operator on $\h$, $Y_1,\dots,Y_n$ are positive contractions on $\h$ summing to $1$ and $v$ is a vector in $\h$.

Theorem \ref{thm1.2} below characterizes those functions on $\Pi^n$ that have a representation of the form \eqref{ctn}. To state this theorem we require a notion based on the following classical result of Pick \cite{Pick}.
 \begin{theorem}\label{pickcond1}
A function  $h$ defined on $\Pi$ belongs to $\Pick$ if and only if the function $A$ defined on $\Pi \times \Pi$ by
  \[
  A(z, w) = \frac{h(z) - \cc{h(w)}}{z - \cc{w}}
  \]
   is positive semidefinite,  that is, for all $n \geq 1, z_1, \dots, z_n \in \Pi, c_1, \dots, c_n \in \C$,
  \[
  \sum A(z_j, z_i)\cc{c_i} c_j \geq 0.
  \]
  \end{theorem}

The following theorem, proved in \cite{ag90}, leads to a generalization of Theorem \ref{pickcond1} to two variables. The {\em Schur class of the polydisc}, denoted by $\schur$, is the set of analytic functions on the polydisc $\D^n$ that are bounded by $1$ in modulus.
 \begin{theorem}\label{schurmodel}
  A function $\ph$  defined on $\D^2$ belongs to $\mathcal S_2$ if and only if there exist positive semidefinite functions $A_1$ and $A_2$ on $\D^2$ such that
  \beq
   1 - \cc{\ph(\mu)}\ph(\la) = (1 - \cc\mu_1\la_1)A_1(\la, \mu)+ (1 - \cc\mu_2\la_2)A_2(\la, \mu).
  \eeq
 \end{theorem}

By way of the transformations
\beq\label{cayleyin}
  z = i\frac{1+\lambda}{1-\lambda},  \quad \la = \frac{z-i}{z+i},
\eeq and
\beq\label{cayleyout}
  h(z) = i\frac{1+\vp(\lambda)}{1-\vp(\lambda)},  \quad \ph(\la)= \frac{h(z)-i}{h(z)+i},
 \eeq
there is a one-to-one correspondence between functions in the  Schur and Pick classes. Under these transformations, Theorem \ref{schurmodel} becomes the following generalization of Pick's theorem to two variables.

 \begin{theorem}\label{pickcond2}
  A function $h$ defined on $\Pi^2$ belongs to $ \Pick_2$ if and only if there exist positive semidefinite functions $A_1$ and $A_2$ on $\Pi^2$ such that
  \[
  h(z) - \cc{h(w)} = (z_1 - \cc{w_1})A_1(z,w) + (z_2 - \cc{w_2})A_2(z,w).
  \]
 \end{theorem}

In the light of Theorems \ref{pickcond1} and \ref{pickcond2} we define the {\em Loewner class} $\p$ to be the set of analytic functions $h$ on $\Pi^n$ with the property that there exist $n$ positive semidefinite functions $A_1,\dots,A_n$ on $\Pi^n$ such that
\beq \label{loewnercond}
h(z) - \overline{h(w)} = \sum_{j=1}^n (z_j-\overline{w_j}) A_j(z,w)
\eeq
for all $z,w\in\Pi^n$.
The Loewner class $\p$ played a key role in \cite{AMY}, which gave a generalization to several variables of Loewner's characterization of the one-variable operator-monotone functions  \cite{lo34}. As the following theorem makes clear, $\p$ also has a fundamental role to play in the understanding of Nevanlinna representations in several variables.
\begin{theorem} \label{thm1.2}
  A function $h$ defined on $\Pi^n$ has a representation of the form \eqref{ctn} if and only if $h \in \p$ and
   \beq\label{strong}
    \liminf_{y\to\infty} y|h(iy,\dots,iy)| < \infty.
   \eeq
 \end{theorem}
In the cases when $n=1$ and $n=2$, Theorems \ref{pickcond1} and \ref{pickcond2} assert that that $\p = \Pick_n$, and so for $n=1$, Theorem \ref{thm1.2} is Nevanlinna's classical Theorem \ref{thm1.0}, and when $n=2$, Theorem \ref{thm1.2} is a straightforward generalization of that result to two variables. When there are more than two variables, it is known that the Loewner class is a proper subset of the Pick class, $\p \neq \Pick_n$ \cite{par70,var74}. Nevertheless, Nevanlinna's result survives as a theorem about the representation of elements of $\p$. Other than the work in  \cite{gkvw}  very little is known about the representation of functions in $\Pick_n$ for three or more variables.

For a function $h$ on $\Pi^n$, we call the formula \eqref{ctn} a {\em Nevanlinna representation of type $1$}. Thus, Theorem \ref{thm1.2} can be rephrased as the assertion that $h$ has a Nevanlinna representation of type $1$ if and only if $h \in \p$ and $h$ satisfies condition \eqref{strong}. Somewhat more complicated representation formulae are needed to generalize Theorem \ref{thm1.1}. We identify three further representation formulae, of increasing generality, and show that every function in $\p$ has a representation of one or more of the four types.

For a function $h$ defined on $\Pi^n$, we refer to a formula
\beq
h(z_1,\dots,z_n) = a + \ip{(A-z_1Y_1-\dots -z_nY_n)\inv v}{v}_\h \qquad \mbox{ for } z_1,\dots,z_n \in\Pi,
\eeq
where $a$ is a constant, $\h$ is a Hilbert space, $A$ is a densely defined self-adjoint operator on $\h$, $Y_1,\dots,Y_n$ are positive contractions on $\h$ summing to $1$ and $v$ is a vector in $\h$, as a \emph{ Nevanlinna representation of type $2$}.
\begin{theorem}\label{type2intro}
A function $h$ defined on $\Pi^n$ has a Nevanlinna representation of type $2$ if and only if $h \in \p$ and
   \beq
    \liminf_{y\to\infty} y\im h(iy,\dots,iy) < \infty.
   \eeq
 \end{theorem}

A {\em Nevanlinna representation of  type $3$} of a function $h$ defined on $\Pi^n$ is of the form
\[
h(z) = a+ \ip{(1-iA)(A - z_Y)^{-1} (1+z_YA)(1-iA)\inv v}{v}   \quad\mbox{ for all } z\in\Pi^n
\]
for some real $a$, some self-adjoint operator $A$ and some vector $v$, where $Y_1, \dots,Y_n$ are operators as in equation \eqref{ctn} above and $z_Y=z_1Y_1+\dots+z_nY_n$.
 \begin{theorem}\label{type3intro}
 A function $h$ defined on $\Pi^n$ has a Nevanlinna representation of type $3$ if and only if $h \in \p$ and
  \[
  \liminf_{y\to\infty} \frac{1}{y}  \im h(iy, \dots, iy)=0.
  \]
 \end{theorem}

Finally, {\em Nevanlinna representations of type $4$} are given by the formula
\beq \label{type4intro}
h(z) = \ip{M(z)v}{v},
\eeq where $M(z)$ is an operator of the form
 \beq \label{Mz}
 \bbm -i&0\\0&1-iA \ebm
   \left( \bbm 1&0 \\ 0 & A \ebm  -
   z_P\bbm 0&0\\0& 1 \ebm \right)\inv
\left(z_P\bbm 1&0\\0& A\ebm +
   \bbm 0&0\\ 0&1 \ebm \right)
   \bbm  -i&0\\0&1-iA \ebm^{-1},
\eeq acting on an orthogonal direct sum of Hilbert spaces $\mathcal N \oplus \M$. In \eqref{type4intro}, $v$ is a vector in $\mathcal N \oplus \M$. In \eqref{Mz}, $A$ is a densely-defined self-adjoint operator acting on $\mathcal M$ and $z_P$ is the operator acting on $\mathcal N \oplus \M$ via the formula
\[
z_P = \sum z_i P_i
\]
where $P_1, \dots, P_n$ are pairwise orthogonal projections acting on $\N \oplus \M$ that sum to $1$.
\begin{theorem}\label{thm2.4}
 Let $h$ be a function defined on $\Pi^n$. Then $h$ has a Nevanlinna representation of type $4$ if and only if $h \in \p$.
\end{theorem}
A weaker, ``generic" version of Theorem \ref{thm2.4} appeared in \cite[Theorem 6.9]{AMY}, where it was used to show that elements in $\p$ are locally operator-monotone.

 It turns out that for $1 \leq k \leq 4$, if $h$ is a function on $\Pi^n$ and $h$ has a Nevanlinna representation of type $k$, then for $k \leq j \leq 4$, $h$ also has a Nevanlinna representation of type $j$. Thus, it is natural to define the \emph{type} of a function in $\p$ to be the smallest $k$ such that $h$ has a Nevalinna representation of type $k$.

For $h \in \p$ the type of $h$ can be characterized in function-theoretic terms through the use of a geometric idea due to Carath\'eodory.  A \emph{carapoint} for a function $\vp$ in the Schur class $\schur$ is a point $\tau \in \T$ such that
\[
\liminf_{\lambda \to \tau} \frac{1 - \abs{\vp(\lambda)}}{1 - \norm{\lambda}_\infty} < \infty,
\]
where
\[
\norm{\lambda}_\infty = \max_{1\leq i \leq n} \abs{\lambda_i}.
\]
Carath\'eodory introduced this notion in one variable in \cite{car29}, along the way to refining earlier results of Julia \cite{ju20}. The following was Carath\'eodory's main result; the notation $\la\nt\tau$ means that $\la$ tends nontangentially to $\tau$.
 \begin{theorem}
 Let $\vp \in \mathcal S_1, \tau \in \T$. If $\tau$ is a carapoint for $\vp$, then $\vp$ is nontangentially differentiable at $\tau$, that is, there exist values $\vp(\tau)$ and $\vp'(\tau)$ such that
 \[
 \lim_{\lambda \nt \tau} \frac{\vp(\lambda) - \vp(\tau) - \vp'(\tau)(\lambda - \tau)}{\lambda - \tau} = 0.
 \]
 In particular, if $\tau$ is a carapoint for $\vp$ then there exists a unique point $\vp(\tau) \in \T$ such that $\vp(\lambda) \to \vp(\tau)$ as $\lambda \nt \tau$.
 \end{theorem}
 In several variables, carapoints have been studied in \cite{abate,jafari,amy10a}. The strong conclusion of nontangential differentiability is lost in several variables; however, at a carapoint $\tau$, there still exists a unimodular nontangential limit $\vp(\tau)$.

 As the point $\chi = (1, \dots, 1)$ is transformed to the point $\infty = (\infty, \dots, \infty)$ by \eqref{cayleyin}, it is natural to say that a function $h \in \pick$ has a carapoint at $\infty$ if the associated Schur function $\ph$, given by the transformation in \eqref{cayleyout}, has a carapoint at $\chi$, and in that case to define $h(\infty)$ by
\beq
h(\infty) = i\frac{1 + \ph(\chi)}{1 - \ph(\chi)}.
\eeq
\black

The connection between carapoints and function types is given in the following theorem.
\begin{theorem}\label{typeoffcn}
 For a function $h \in \p$,
\begin{enumerate}
\item $h$ is of type $1$ if and only if $\infty$ is a carapoint of $h$ and $h(\infty) = 0$;
 \item $h$ is of type $2$ if and only if $\infty$ is a carapoint of $h$ and $h(\infty)\in \R\setminus\{0\}$;
 \item $h$ is of type $3$ if and only if $\infty$ is not a carapoint of $h$;
 \item $h$ is of type $4$ if and only if $\infty$ is a carapoint of $h$ and $h(\infty) = \infty$.
\end{enumerate}
\end{theorem}

The paper is structured as follows. As is clear from the formulae used to define the various Nevanlinna representations,  Nevanlinna representations are generalizations of the resolvent of a self-adjoint operator.  These \emph{structured resolvents}, studied in Sections \ref{structured} and \ref{matricial}, are analytic operator-valued functions on the polyhalfplane $\Pi^n$ with non-negative imaginary part, fully analogous to the familiar resolvent operator.  There are also {\em structured resolvent identities} for them, studied in Section \ref{resolvident} of the paper.

In modern texts Nevanlinna's representation is derived from the Herglotz Representation with the aid of the Cayley transform \cite{lax02,don74}. In Section \ref{type4} we introduce the $n$-variable \emph{strong Herglotz class} and then prove Theorem \ref{type4intro} by applying the Cayley transform to Theorem 1.8 of \cite{ag90}.

In Section \ref{type321} we derive the Nevanlinna representations of type $3, 2$, and $1$, we show how they arise naturally  from the underlying Hilbert space geometry and we prove slight strengthenings of Theorems \ref{thm1.2}, \ref{type2intro} and \ref{type3intro}.  In Section \ref{asymptotic} we give function-theoretic conditions for a function $h \in \p$ to possess a representation of a given type.

In Section \ref{caraInfty} we introduce the notion of carapoints for functions in the Pick class and in Section \ref{caraTypes} we establish the criteria in Theorem \ref{typeoffcn} for the type of a function using the language of carapoints.

In Section \ref{growth} we give the growth estimates  for functions in $\p$ that flow from our analysis of structured resolvents, and in Section \ref{resolvident} we present resolvent identities for structured resolvents.

Results related to ours from a system-theoretic perspective have been obtained in ongoing work of J. A. Ball and D. Kalyuzhnyi-Verbovetzkyi \cite{BKV1,BKV2}.  See also \cite{BS}, where Krein space methods are applied to similar problems.

\section{Structured resolvents of operators}\label{structured}
The resolvent operator $(A-z)\inv$ of a densely defined self-adjoint operator $A$ on a Hilbert space plays a prominent role in spectral theory.  It has the following properties.  
\begin{enumerate}
\item It is an analytic bounded operator-valued function of $z$ in the upper halfplane $\Pi$;
\item it satisfies the growth estimate $\|(A-z)\inv\| \leq 1/ \im z$ for $z\in\Pi$;
\item $(A-z)\inv$ has non-negative imaginary part for all $z\in\Pi$;
\item it satisfies the ``resolvent identity".    
\end{enumerate}
Here we are interested in several-variable analogs of the resolvent.  These will again be operator-valued analytic functions with non-negative imaginary part, but now on the polyhalfplane $\Pi^n$.  Because of the additional complexities in several variables we encounter three different types of resolvent; all of them have the four listed properties, with very slight modifications, and therefore deserve the name {\em structured resolvent}.

For any Hilbert space $\h$, a {\em positive decomposition} of $\h$ will mean an $n$-tuple $Y=(Y_1,\dots,Y_n)$ of positive contractions on $\h$ that sum to the identity operator. For any $z=(z_1,\dots,z_n)\in\C^n$  and any $n$-tuple $T=(T_1,\dots,T_n)$ of bounded operators we denote by $z_T$ the operator $\sum_j z_j T_j $.  Here each $T_j$ is a bounded operator from $\h_1$ to $\h_2$, for some Hilbert spaces $\h_1, \ \h_2$, so that $z_T$ is also a bounded operator from $\h_1$ to $\h_2$.
\begin{definition} \label{resolv2}
Let  $A$ be a closed densely defined self-adjoint operator on a Hilbert space $\h$ and let $Y$ be a positive decomposition of $\h$.  The {\em structured resolvent of $A$ of type 2} corresponding to $Y$ is the operator-valued function 
\[
 z\mapsto (A-z_Y)\inv : \Pi^n \to  \L(\h).
\]
\end{definition}
The following observation is essentially \cite[Lemma 6.25]{AMY}.
\begin{proposition}\label{AzY}
For $A$ and $Y$ as in Definition {\rm \ref{resolv2}} the structured resolvent $(A-z_Y)\inv$ is well defined on $\Pi^n$ and satisfies, for all $z\in\Pi^n$,
\beq\label{boundAzYinv}
\|(A-z_Y)\inv\| \leq \frac{1}{\min_j \im z_j}.
\eeq
Moreover
\begin{align} \label{imPos}
\im \left((A-z_Y)\inv\right) &= (A-z_Y^*)\inv \left(\im z_Y\right) (A-z_Y)\inv \\
	&= (A-z_Y)\inv \left(\im z_Y\right) (A-z_Y^*)\inv \nn\\
	&\geq 0. \nn
\end{align}
\end{proposition}
The range of the bounded operator $(A-z_Y)^{-1}$ is of course $\mathcal D(A)$, the domain of $A$. 
\begin{proof}
For any vector $\xi$ in the domain of $A$,
\begin{align*}
\|(A-z_Y)\xi\| \ \|\xi\| &\geq |\ip{(A-z_Y)\xi}{\xi}| \\
	&\geq |\im \ip{(A-z_Y)\xi}{\xi}| \\
	&= \ip{(\im z_Y) \xi}{\xi} \\
	&= \sum_j (\im z_j)\ip{Y_j\xi}{\xi} \\
	&\geq (\min_j  \im z_j) \ip{\sum_j Y_j\xi}{\xi}\\
	&= (\min_j  \im z_j) \|\xi\|^2.
\end{align*}
Thus $A-z_Y$ has lower bound $\min_j\im z_j > 0$, and so has a bounded left inverse.  A similar argument with $z$ replaced by $\bar z$ shows that $(A-z_Y)^*$ also has a bounded left inverse, and so $A-z_Y$ has a bounded inverse and the inequality \eqref{boundAzYinv} holds.

The identities \eqref{imPos} are easy.
\end{proof}

Resolvents of type 2 are the simplest several-variable analogues of the familiar one-variable resolvent but they are not sufficient for the analysis of the several-variable Pick class.  To this end we introduce two further generalizations.  Let us first recall some basic facts about closed unbounded operators.
\begin{lemma}\label{basicUnbdd}
Let $T$ be a closed densely defined operator on a Hilbert space $\h$, with domain $\mathcal D(T)$.  The operator $1+T^*T$ is a bijection from $\mathcal D(T^*T)$ to $\h$, and the operators
\[
B\df (1+T^*T)^{-1}, \qquad C\df T(1+T^*T)^{-1}
\]
are everywhere defined and contractive on $\h$.  Moreover $B$ is self-adjoint and positive, and 
$ \ran C \subset \mathcal D(T^*)$.
\end{lemma}
\begin{proof}
All these statements are proved in \cite[Sections 118, 119]{RN}, although the final statement about $\ran C$ is not explicitly stated.  We must show that for all $v\in\h$ there exists $y\in\h$ such that, for all $h\in\h$,
\[
\ip{Th}{Cv}=\ip{h}{y}.
\]
It is straightforward to check that this relation holds for $y=v-Bv$, and so $\ran C\subset \mathcal D(T^*)$.
\end{proof}

\begin{definition} \label{resolv3}
Let  $A$ be a closed densely defined self-adjoint operator on a Hilbert space $\h$ and let $Y$ be a positive decomposition of $\h$.  The {\em structured resolvent of $A$ of type 3} corresponding to $Y$ is the operator-valued function 
$M:\Pi^n\to  \L(\h)$ given by
\beq\label{defM3}
M(z) = (1-iA)(A - z_Y)^{-1} (1+z_Y A)(1-iA)^{-1}.
\eeq
\end{definition}
We denote the $\ell_1$ norm on $\C^n$ by $\|\cdot\|_1$.  Note that $\|z_Y\| \leq \|z\|_1$ for all $z\in\C^n$ and all positive decompositions $Y$.
\begin{proposition}\label{existResolv3}
For $A$ and $Y$ as in Definition {\rm \ref{resolv3}} the structured resolvent $M(z)$ of type $3$ given by equation \eqref{defM3} is well defined as a bounded operator on $\h$  for all $z\in\Pi^n$ and satisfies
\beq\label{boundM3}
\|M(z)\| \leq  (1+2\|z\|_1) \left(1+ \frac{1+\|z\|_1}{\min_j \im z_j}\right).
\eeq
\end{proposition}

\begin{proof}
Since
\[
1+z_YA = 1-iz_Y + iz_Y(1-iA) : \mathcal D(A) \to \h
\]
and $(1-iA)\inv$ is a contraction on all of $\h$, with range $\mathcal D(A)$, the operator $(1+z_YA)(1-iA)\inv$ is well defined as an operator on $\h$ and
\begin{align}\label{3&4}
\|(1+z_YA)(1-iA)\inv\| &= \| (1-iz_Y)(1-iA)\inv + iz_Y\| \nn \\
	&\leq \|1-iz_Y\| + \|z_Y\| \nn\\
	&\leq 1+2\|z_Y\| \nn\\
	&\leq 1+2\|z\|_1.
\end{align}
Similarly $(1-iA)(A-z_Y)\inv$ is well defined on $\h$, and since
\[
i(A-z_Y)= -(1-iA)+(1-iz_Y): \mathcal D(A) \to \h
\]
we have
\[
i=-(1-iA)(A-z_Y)\inv + (1-iz_Y)(A-z_Y)\inv:\h\to\h.
\]
Thus, by virtue of the bound \eqref{boundAzYinv},
\begin{align}\label{1&2}
\|(1-iA)(A-z_Y)\inv\| &= \|i - (1-iz_Y)(A-z_Y)\inv\| \nn \\
	&\leq 1+ \|1-iz_Y\| \  \|(A-z_Y)\inv\|  \nn \\
	&\leq 1+ \frac{1+\|z\|_1}{\min_j \im z_j}.
\end{align}
On combining the estimates \eqref{1&2} and \eqref{3&4} we obtain the bound \eqref{boundM3}.
\end{proof}
The following alternative formula for the structured resolvent of type 3, valid on the dense subspace $\mathcal D(A)$ of $\h$, allows us to show that $\im M(z) \geq 0$.
\begin{proposition}\label{2ndResolv3}
For $A$ and $Y$ as in Definition {\rm \ref{resolv3}} and $z\in\Pi^n$
\begin{align} 
M(z)|\mathcal D(A) &= (1-iA)\left\{(A-z_Y)\inv-A(1+A^2)\inv\right\}(1+iA) \label{Resolv3.2}\\
	&= (1-iA)(A-z_Y)\inv(1+iA) - A:  \mathcal D(A) \to \h. \label{3rdtype3}
\end{align}
Moreover, for every $v\in\mathcal D(A)$,
\beq\label{ImMpos}
\im \ip{M(z)v}{v} = \ip {(1-iA)(A-z_Y^*)\inv (\im z_Y) (A-z_Y)\inv(1+iA)v}{v} \geq 0.
\eeq
\end{proposition}
\begin{proof}
By Lemma \ref{basicUnbdd} the operator $A(1+A^2)\inv$ is contractive on $\h$ and has range contained in $\mathcal D(A)$.   On $\mathcal D(A^2)$ we have the identity
\[
1+z_YA= 1+A^2- (A-z_Y)A.
\]
Since $(1+A^2)\inv$ maps $\h$ into $\mathcal D(A^2)$ we have
\[
(1+z_Y A)(1+A^2)\inv = 1-(A-z_Y)A(1+A^2)\inv : \h \to \h,
\]
and therefore
\beq\label{3f}
(A-z_Y)\inv(1+z_Y A)(1+A^2)\inv = (A-z_Y)\inv-A(1+A^2)\inv : \h \to \mathcal D(A).
\eeq
Clearly
\[
 (1+A^2)\inv (1+iA)=(1-iA)\inv  \quad \mbox{ on }\mathcal D(A)
\]
and so, on multiplying equation \eqref{3f} fore and aft by $1\pm iA$, we deduce that, as operators from $\mathcal D(A)$ to $\h$,
\begin{align*}
M(z) | \mathcal D(A)&= (1-iA)(A-z_Y)\inv (1+z_YA)(1-iA)\inv  \nn\\
	&= (1-iA)(A-z_Y)\inv (1+z_YA)(1+A^2)\inv(1+iA) \nn\\
	&=(1-iA)\left\{ (A-z_Y)\inv-A(1+A^2)\inv \right\}(1+iA).  
\end{align*}
This establishes equation \eqref{Resolv3.2}.

The expression \eqref{3rdtype3} follows from equation \eqref{Resolv3.2} since
\[
(1-iA)A(1+A^2)\inv (1+iA) =A  \quad \mbox{ on } \mathcal D(A).
\]

By equation \eqref{3rdtype3} we have, for any $z\in\Pi^n$ and $v\in \mathcal D(A)$,
\begin{align*}
\im \ip{M(z)v}{v} &= \im\ip{(1-iA)(A-z_Y)\inv(1+iA)v}{v} -\im \ip{Av}{v}\\
	&= \im\ip{(A-z_Y)\inv(1+iA)v}{(1+iA)v}
\end{align*}
and hence, by equation \eqref{imPos},
\[
\im \ip{M(z)v}{v} = \ip{(A-z_Y^*)\inv (\im z_Y) (A-z_Y)\inv(1+iA)v}{(1+iA)v},
\]
and so equation \eqref{ImMpos} holds.
\end{proof}
\begin{corollary} \label{Mpick}
For $A$ and $Y$ as in Definition {\rm \ref{resolv3}} the structured resolvent $M(z)$ given by equation \eqref{defM3} satisfies $\im M(z) \geq 0$ for all $z\in\Pi^n$.
\end{corollary}
For, by Propositions \ref{existResolv3} and \ref{2ndResolv3},  $M(z)$ is a bounded operator on $\h$ and $\im \ip{M(z)v}{v} \geq 0$ for $v\in\mathcal D(A)$.  The conclusion follows by density of $\mathcal D(A)$ and continuity.

In the case of bounded $A$ there is yet another expression for the structured resolvent of type $3$.
\begin{proposition}
If $A $ is a {\em bounded} self-adjoint operator on $\h$  and $Y$ is a positive decomposition of $\h$ then, for $z\in\Pi^n$,
\beq\label{Resolv3.3}
M(z)=(1+iA)\inv(1+Az_Y)(A-z_Y)\inv(1+iA)
\eeq
\end{proposition}
\begin{proof}
Since $A$ is bounded it is defined on all of $\h$.  We have
\[
1+Az_Y = 1+ A^2 - A(A-z_Y)
\]
and hence
\[
(1+Az_Y)(A-z_Y)\inv = (1+ A^2)(A-z_Y)\inv -A.
\]
Thus
\begin{align*}
(1+iA)\inv(1+Az_Y)(A-z_Y)\inv(1+iA) &=  (1-iA)(A-z_Y)\inv(1+iA) -A \\
	&=M(z)
\end{align*}
by equation \eqref{3rdtype3}.
\end{proof}

\begin{remark} \rm
In the case of unbounded $A$ the expression \eqref{Resolv3.3} for $M(z)$ is valid wherever it is defined, but it is 
not to be expected that this will be a dense subspace of $\h$ in general.
\end{remark}

Here are two examples of structured resolvents of type 3, one on $\C^2$ and one on an infinite-dimensional space.
\begin{example} \label{type3ex1} \rm
Let 
\[
\h=\C^2,\quad A=\bbm 1&0\\0&-1 \ebm, \quad Y_1=\half \bbm 1&1\\1&1 \ebm, \quad Y_2 =1-Y_1, \quad Y=(Y_1,Y_2).
\]
Then 
\begin{align*}
M(z) &= (1-iA)(A-z_Y)^{-1}(1+z_Y A)(1-iA)^{-1} \\
	&= \frac{1}{1-z_1z_2}\bbm (1+z_1)(1+z_2) & -i(z_1-z_2) \\ i(z_1-z_2) & -(1-z_1)(1-z_2) \ebm.
\end{align*}
\end{example}
\begin{example}\label{type3ex2} \rm
Let $\h=L^2(\R)$,  let $A$ be the operation of multiplication by the independent variable $t$ and let $Y=(P,Q)$ where $P, Q$ are the orthogonal projection operators onto the subspaces of even and odd functions respectively in $L^2$.  Thus
\[
Pf(t)=\half\left\{f(t)+f(-t)\right\}, \qquad Qf(t)=\half\left\{f(t)-f(-t)\right\}.
\]
Let $Y'=(Q,P)$.  Note that 
\[
PA=AQ, \qquad QA=AP
\]
and hence
\[
z_YA=Az_{Y'}, \qquad z_{Y'}A=Az_Y, \qquad z_Yz_{Y'}=z_1z_2=z_{Y'}z_Y.
\]
It follows that $z_Y$ and $z_{Y'}$ commute with $A^2$, and it may be checked that
\[
(A-z_Y)\inv=(A^2-z_1z_2)\inv(z_{Y'}+A)=(z_{Y'}+A)(A^2-z_1z_2)\inv
\]
and hence
\[
(A-z_Y)\inv(1+z_YA)= (A^2-z_1z_2)\inv\left( (1+A^2)z_{Y'} + (1+z_1z_2)A\right).
\]
A straightforward calculation now shows that the structured resolvent $M(z)$ of $A$ corresponding to $Y$ is given by
\[
(M(z)f)(t)= \frac{ \left( \half (z_1+z_2)(1+t^2)+(1+z_1z_2)t\right)f(t) + \half(z_2-z_1)(1-it)^2f(-t)}{t^2-z_1z_2}
\]
for all $z\in\Pi^2,\ f\in L^2(\R)$ and $t\in\R$.
In particular, we note for future use that if $f$ is an even function,
\beq\label{fcnoftype3}
(M(z)f)(t)= \frac{t(1+z_1z_2)+(1-it)(itz_1+z_2)}{t^2-z_1z_2}f(t).
\eeq
\end{example}

\section{The matricial resolvent} \label{matricial}
The third and last form of structured resolvent that we consider has a $2\times 2$ matricial form.  As will become clear, this extra complication is needed for the description of the most general type of  function in the several-variable Loewner class.

By an {\em orthogonal decomposition} of a Hilbert space $\h$ we shall mean an $n$-tuple $P=(P_1,\dots,P_n)$ of orthogonal projection operators with pairwise orthogonal ranges such that $\sum_{j=1}^n P_j$ is the identity operator. 
\begin{proposition}\label{2x2resolv}
Let $\mathcal H$ be the orthogonal direct sum of Hilbert spaces $\mathcal{ N, M}$, let $A$ be a densely defined self-adjoint operator on $\mathcal M$ with domain $\mathcal{D}(A) $ and let $P$ be an orthogonal decomposition of $\mathcal H$.  For every $z\in\Pi^n$ the operator on $\mathcal H$ given with respect to the decomposition $\mathcal {N \oplus M}$ by the matricial formula 
\begin{align}\label{defM}
M(z)&=\bbm -i&0\\0&1-iA \ebm
   \left( \bbm 1&0 \\ 0 & A \ebm  -
   z_P\bbm 0&0\\0& 1 \ebm \right)\inv 
\left(z_P\bbm 1&0\\0& A\ebm +
   \bbm 0&0\\ 0&1 \ebm \right)
   \bbm  -i&0\\0&1-iA \ebm^{-1}
\end{align}
is a bounded operator defined on all of $\h$, and
\beq\label{estNormM}
\|M(z)\| \leq (1+\sqrt{10} \|z\|_1)\left(1 + \frac{1+\sqrt{2}\|z\|_1}{\min_j \im z_j}\right)
\eeq
\end{proposition}

\begin{proof}
Let $z\in\Pi^n$.
Let  the projection $P_j$ have operator matrix
\beq\label{expP}
P_j = \begin{bmatrix} X_j & B_j \\ B_j^* & Y_j \end{bmatrix}
\eeq
with respect to the decomposition $\mathcal{H = N \oplus M}$.  Then 
\[
X=(X_1,\dots,X_n), \quad Y=(Y_1, \dots, Y_n)
\]
are positive decompositions of $\N, \ \M$ respectively, and
\[
B=(B_1, \dots, B_n), \quad B^*=(B_1^*, \dots, B_n^*)
\]
are $n$-tuples of contractions summing to $0$, from $\M$ to $\N$ and from $\N$ to $\M$ respectively.  Since the $B_j$ are contractions we have
\[
\|z_B\| \leq \|z\|_1.
\]
For any $z\in\C^n$,
\beq\label{formzP}
z_P = \begin{bmatrix} z_X & z_B \\ z_{B^*} & z_Y \end{bmatrix}.
\eeq

Consider the third and fourth factors in the product on the right hand side of equation \eqref{defM}; the product of these two factors is well defined as an operator on $\h$ since $(1-iA)\inv$ maps $\M$ to $\mathcal D(A)$.  It is even a bounded operator, since, by virtue of equation \eqref{formzP},
\begin{align}\label{last2}
 \left(z_P\bbm 1&0\\0& A\ebm +
   \bbm 0&0\\ 0&1 \ebm \right)
   \bbm  -i&0\\0&1-iA \ebm^{-1}  &=
	\bbm iz_X & z_BA(1-iA)\inv \\ iz_{B^*} & (1+z_YA)(1-iA)\inv \ebm.
\end{align}
Since
\[
\|A(1-iA)\inv\|= \|i\left(1-(1-iA)\inv\right)\| \leq 2
\]
we can immediately see that the operator \eqref{last2} is bounded. We can get an estimate of the norm of the operator matrix \eqref{last2} if we replace each of the four operator entries by an upper bound for its norm.  We find that 
\begin{align}\label{normlast2}
\left\|\left(z_P\bbm 1&0\\0& A\ebm +
   \bbm 0&0\\ 0&1 \ebm \right)
   \bbm  -i&0\\0&1-iA \ebm^{-1}\right\| &\leq \left\| \bbm \|z\|_1 & 2\|z\|_1 \\ \|z\|_1 & 1+2\|z\|_1 \ebm \right\| \nn \\
	&\leq 1+\|z\|_1\left\|\bbm 1 & 2\\1& 2\ebm \right\|\nn  \\
	&= 1+\sqrt{10} \|z\|_1.
\end{align}

Now consider the second factor in the definition \eqref{defM} of $M(z)$.  We find that
\begin{align}\label{needThis}
\left( \bbm 1&0 \\ 0 & A \ebm  -
   z_P\bbm 0&0\\0& 1 \ebm\right)\inv  &= \bbm 1 & -z_B \\ 0 & A-z_Y \ebm\inv \nn \\
	&= \bbm 1& z_B(A-z_Y)\inv \\ 0& (A-z_Y)\inv \ebm,
\end{align}
which maps $\h$ into $\N\oplus \mathcal D(A)$.  Hence the product of the first two factors in the product on the right hand side of equation \eqref{defM} is
\beq\label{first2}
\bbm -i&0\\0&1-iA \ebm
   \left( \bbm 1&0 \\ 0 & A \ebm  -
   z_P\bbm 0&0\\0& 1 \ebm \right)\inv = \bbm -i & -iz_B(A-z_Y)\inv \\ 0& (1-iA)(A-z_Y)\inv \ebm.
\eeq
Since
\begin{align*}
\|(1-iA)(A-z_Y)\inv\| &= \|(1-iz_Y)(A-z_Y)\inv -i\| \\
	&\leq 1+ \|1-iz_Y\| \, \|(A-z_Y)\inv\| \\
	&\leq 1+\frac{1+\|z\|_1}{\min_j \im z_j}
\end{align*}
we deduce from equation \eqref{first2} that 
\begin{align}\label{normfirst2}
\left\| \bbm -i&0\\0&1-iA \ebm
   \left( \bbm 1&0 \\ 0 & A \ebm  -
   z_P\bbm 0&0\\0& 1 \ebm \right)\inv  \right\| &\leq  \left\|\bbm 1& \|z\|_1 \ \|(A-z_Y)\inv \| \\ 0 & 1+(1+\|z\|_1) \|(A-z_Y)\|\inv\ebm \right\| \nn \\
	&  \leq 1 +\left\|\bbm 0&\|z\|_1 \\ 0& 1+\|z\|_1 \ebm \bbm 0&0\\0& \|(A-z_Y)\inv\| \ebm \right\| \nn \\
	&\leq  1 + \frac{1+\sqrt{2}\|z\|_1}{\min_j \im z_j}.
\end{align}
On combining the estimates \eqref{normfirst2} and \eqref{normlast2} we obtain the bound \eqref{estNormM} for $\|M(z)\|$.
\end{proof}
\begin{remark}\label{3.2} \rm
On multiplying together the expressions \eqref{first2} and \eqref{last2} we obtain the formula
\[
M(z)= \bbm z_X+z_B(A-z_Y)\inv z_{B^*}  & -iz_B(A-z_Y)\inv(1+iA) \\
	i(1-iA)(A-z_Y)\inv z_{B^*} & (1-iA)(A-z_Y)\inv(1+z_YA)(1-iA)\inv \ebm.
\]
Notice in particular that the $(2,2)$ entry (that is, the compression of $M(z)$ to $\M$) is the structured resolvent of $A$ of type $3$ corresponding to $Y$, the compression of $P$ to $\M$, as in equation \eqref{defM3}.
\end{remark}
\begin{definition}\label{defMatResolv}
Let $\mathcal H$ be the orthogonal direct sum of Hilbert spaces $\mathcal{ N, M}$, let $A$ be a densely defined self-adjoint operator on $\mathcal M$ with domain $\mathcal{D}(A) $ and let $P$ be an orthogonal decomposition of $\mathcal H$.  The {\em structured resolvent of $A$ of type $4$} corresponding to $P$ is the operator-valued function $M: \Pi^n \to \L(\h)$ given by equation \eqref{defM}.
\end{definition}
We shall also refer to $M(z)$ as the {\em matricial resolvent of $A$ with respect to $P$}. The important property that $\im M(z)\geq 0$ is not at once apparent from the formula \eqref{defM};  as with structured resolvents of type $3$,   there are alternative formulae from which this property is more easily shown.  Once again the alternatives suffer the minor drawback that they give $M(z)$ only on a dense subspace of $\h$.
\begin{proposition}\label{alternM}
With the notation of Definition {\rm \ref{defMatResolv}},  as operators on $ \N\oplus\mathcal D(A)$,
\begin{align}\label{2ndM}
M(z)&=\bbm -i&0\\0&1-iA \ebm \left(\bbm 1&0\\0& A(1+A^2)\inv\ebm z_P + \bbm 0&0\\0&(1+A^2)\inv \ebm\right)
  \times \nn \\		&\hspace*{3cm}
\left(\bbm 1&0\\0&A \ebm - \bbm 0&0\\0&1\ebm z_P \right)\inv \bbm i&0\\0&1+iA\ebm\\
   &= \bbm -i&0\\0&1-iA \ebm  \left(\bbm 1&0\\0&0\ebm z_P + \bbm 0&0\\0&1 \ebm\right)
	\left( \bbm 1&0\\0&A \ebm - \bbm 0&0\\0&1 \ebm z_P\right)\inv  \bbm i&0\\0&1+iA\ebm \nn \\
	& \hspace*{3cm}    - \bbm0&0\\0&A \ebm  \label{simpler} \\
  &= \bbm -i&0\\0&1-iA \ebm  \left( \bbm 1&0\\0&A \ebm - z_P \bbm 0&0\\0&1 \ebm \right)\inv \left(z_P\bbm 1&0\\0&0\ebm  + \bbm 0&0\\0&1 \ebm\right)
	  \bbm i&0\\0&1+iA\ebm \nn \\
	& \hspace*{3cm}    - \bbm0&0\\0&A \ebm  \label{simpler2} 
\end{align}
for all $z\in\Pi^n$. 
 Moreover,  for all $z,\ w\in\Pi^n$,
\begin{align}\label{imM}
 M(z) - M(w)^*&= \bbm -i & 0\\ 0&1-iA \ebm \left(\bbm 1&0\\0&A \ebm - w_P^* \bbm 0&0\\0&1 \ebm\right)^{-1}  \times \nn\\ 
	&\hspace*{2cm} ( z_P- w_P^*) \left(\bbm 1&0\\0&A \ebm -  \bbm 0&0\\0&1 \ebm z_P\right)^{-1} \bbm i&0\\0&1+iA \ebm 
\end{align}
on $\N\oplus \mathcal D(A)$.
\end{proposition}
\begin{proof}
By Lemma \ref{basicUnbdd} the operators $(1+A^2)^{-1}$ and 
\[
 C\df \im (1-iA)^{-1} = A(1+A^2)^{-1}
\]
are self-adjoint contractions defined on all of $\M$.   Furthermore, 
\[
\ran (1+A^2)^{-1} = \mathcal D(A^2), \qquad \ran C \subset \mathcal D(A).
\]
 
We claim that, as operators on $\N\oplus\mathcal D(A)$,
\begin{align}\label{swap}
\left( \bbm 1&0\\ 0&A \ebm \right . & \left . -z_P \bbm 0&0\\0&1 \ebm \right)\inv \left(z_P\bbm 1&0\\0&A\ebm +\bbm 0&0\\0&1 \ebm \right) = \nn \\
	&\left(\bbm 1&0\\0& C\ebm z_P + \bbm 0&0\\0&(1+A^2)\inv \ebm\right)\left( \bbm 1&0\\0&C\ebm - \bbm 0&0\\0&(1+A^2)\inv \ebm z_P\right)\inv.
\end{align}
We have
\begin{align*}
\left(z_P\bbm 1&0\\0&A\ebm +\bbm 0&0\\0&1 \ebm \right)&\left( \bbm 1&0\\0&C\ebm - \bbm 0&0\\0&(1+A^2)\inv \ebm z_P\right)  \\
	& \hspace{-2cm} =\bbm 0&0\\0&C\ebm + z_P \bbm 1&0\\ 0&AC \ebm - \bbm 0&0\\0&(1+A^2)\inv \ebm z_P - z_P\bbm 0&0\\0&C \ebm z_P \\
	& \hspace{-2cm} =\bbm 0&0\\0&C\ebm + z_P \left(\bbm 1&0\\ 0&AC \ebm-1\right)+\left(1 - \bbm 0&0\\0&(1+A^2)\inv \ebm\right) z_P - z_P\bbm 0&0\\0&C \ebm z_P \\
	& \hspace{-2cm} =\bbm 0&0\\0&C\ebm - z_P \bbm 0&0\\0&(1+A^2)\inv \ebm + \bbm 1&0\\0& AC \ebm z_P -z_P\bbm 0&0\\0&C\ebm z_P \\
	& \hspace{-2cm} =\left(\bbm 1&0\\0&A\ebm - z_P \bbm 0&0\\0&1\ebm\right)
	\left(\bbm 1&0\\0&C \ebm z_P + \bbm 0&0 \\0& (1+A^2)\inv \ebm \right).
\end{align*}
This is an identity between operators on $\h$, in both cases a composition $\h\to\N\oplus\mathcal D(A)\to \h$, and moreover the first factor on the left hand side and the second factor on the right hand side are invertible, from $\N\oplus\mathcal D(A)$ to $\h$ and from $\h$ to $\N\oplus\mathcal D(A)$ respectively.
We may pre- and post-multiply appropriately to obtain equation \eqref{swap}, but note that the equation is then only valid as an identity between operators on $\N\oplus\mathcal D(A)$.

On combining equations \eqref{defM} and \eqref{swap} we deduce that
\begin{align*}
M(z)&=\bbm -i&0\\0&1-iA \ebm \left(\bbm 1&0\\0& C\ebm z_P + \bbm 0&0\\0&(1+A^2)\inv \ebm\right) \times \\	
	&\hspace*{3cm}\left( \bbm 1&0\\0&C\ebm - \bbm 0&0\\0&(1+A^2)\inv \ebm z_P\right)\inv \bbm -i&0\\0&1-iA\ebm\inv.
\end{align*}
Since
\[
 \bbm -i&0\\0&1-iA\ebm\inv = \bbm 1&0\\0&1+A^2 \ebm\inv \bbm i&0\\0&1+iA \ebm
\]
and 
\[
\bbm 1&0\\0&1+A^2 \ebm \left( \bbm 1&0\\0&C\ebm - \bbm 0&0\\0&(1+A^2)\inv \ebm z_P\right) =
	\bbm 1&0\\0&A \ebm - \bbm 0&0\\0&1\ebm z_P,
\]
we deduce further that
\begin{align} \label{halfway}
M(z)&=\bbm -i&0\\0&1-iA \ebm \left(\bbm 1&0\\0& C\ebm z_P + \bbm 0&0\\0&(1+A^2)\inv \ebm\right)
  \times \nn \\		&\hspace*{3cm}
\left(\bbm 1&0\\0&A \ebm - \bbm 0&0\\0&1\ebm z_P \right)\inv \bbm i&0\\0&1+iA\ebm,
\end{align}
which proves equation \eqref{2ndM}.
It is straightforward to verify that
\begin{align}\label{stfwd}
\left(\bbm 1&0\\0& C\ebm z_P + \bbm 0&0\\0&(1+A^2)\inv \ebm\right)
	&\left(\bbm 1&0\\0&A \ebm - \bbm 0&0\\0&1\ebm z_P \right)\inv  \\
	&\hspace*{-4cm}= 
	\left(\bbm 1&0\\0&0\ebm z_P + \bbm 0&0\\0&1 \ebm\right)
	\left( \bbm 1&0\\0&A \ebm - \bbm 0&0\\0&1 \ebm z_P\right)\inv - 
	  \bbm0&0\\0&A(1+A^2)\inv\ebm.
\end{align}
Clearly
\[
\bbm -i&0\\0&1-iA\ebm \bbm0&0\\0&A(1+A^2)\inv\ebm \bbm i&0\\0&1+iA\ebm = \bbm 0&0\\0&A \ebm,
\]
and so on suitably pre- and post-multiplying equation \eqref{stfwd}, we obtain equation \eqref{simpler}.

To prove equation \eqref{simpler2}, check first that
\begin{align*}
\left( \bbm 1&0\\0&A \ebm - z_P \bbm 0&0\\0&1\ebm  \right)\left(\bbm 1&0\\0&0 \ebm z_P + \bbm 0&0\\0&1\ebm \right)&=  \\ 
\left(z_P\bbm 1&0\\0&0\ebm+\bbm 0&0\\0&1\ebm \right)&\left(\bbm 1&0\\0&A\ebm - \bbm 0&0\\0&1 \ebm z_P\right)
\end{align*}
as operators on $\N\oplus\mathcal{D}(A)$.
It follows that 
\begin{align*}
\left(\bbm 1&0\\0&0 \ebm z_P + \bbm 0&0\\0&1\ebm \right)\left(\bbm 1&0\\0&A\ebm - \bbm 0&0\\0&1 \ebm z_P\right)\inv &= \\
\left( \bbm 1&0\\0&A \ebm - z_P \bbm 0&0\\0&1\ebm  \right)\inv &\left(z_P\bbm 1&0\\0&0\ebm+\bbm 0&0\\0&1\ebm \right)
\end{align*}
as operators from $\h$ to $\N\oplus \mathcal{D}(A)$.  On combining this relation with equation \eqref{simpler} we derive the expression \eqref{simpler2} for $M(z)|\N\oplus\mathcal{D}(A)$.

We now  derive the identity \eqref{imM}.  Let 
\[
D=\bbm i&0\\0&1+iA \ebm
\]
and consider $z,\ w\in\Pi^n$.  By equation \eqref{2ndM}
\beq\label{MandW}
M(z)=D^* W(z) D
\eeq
 on $\N\oplus \mathcal D(A)$, where 
\begin{align}\label{defW}
W(z) &= R(z)S(z)\inv  - 
 \begin{bmatrix}
  0 & 0 \\ 0 & A(1 + A^2)^{-1}
 \end{bmatrix}
\end{align}
and 
\[
R(z)= \begin{bmatrix}
  1 & 0 \\ 0 & 0
 \end{bmatrix} z_P+
\begin{bmatrix}
     0 & 0 \\ 0 & 1
    \end{bmatrix},  \quad S(z) = \begin{bmatrix}
  1 & 0 \\ 0 & A
 \end{bmatrix} - 
 \begin{bmatrix}
  0 & 0 \\ 0 & 1
 \end{bmatrix} z_P.
\]
We have seen that  $S(z)$ is invertible for any $z\in \Pi^n$, so that $W(z)$ is a bounded operator on $\h$.
Clearly
\begin{align*}
M(z)-M(w)^*&=D^*\left(R(z)S(z)\inv - S(w)^{*-1}R(w)^*\right)D\\
	&= D^*S(w)^{*{-1}}\left( S(w)^*R(z) - R(w)^*S(z)\right) S(z)^{-1}D.
\end{align*}
Here
\begin{align*}
S(w)^*R(z) - R(w)^*S(z) &= \begin{bmatrix} 1&0\\0&0\end{bmatrix} z_P +\begin{bmatrix}0&0\\0&A\end{bmatrix}- w_P^*\begin{bmatrix} 0&0\\ 0&1\end{bmatrix} -\\
	&\hspace{2cm}	\left( w_P^*\begin{bmatrix}1&0\\0&0\end{bmatrix} + \begin{bmatrix} 0&0\\0&A\end{bmatrix} - \begin{bmatrix} 0&0\\0&1\end{bmatrix}z_P \right)\\
	&=z_P-w_P^*.
\end{align*}
Hence
\[
M(z)-M(w)^*= D^*S(w)^{*{-1}}(z_P-w_P^*) S(z)^{-1}D,
\]
which is  equation \eqref{imM}.  
\end{proof}
The next result shows that the matricial resolvent belongs not just to the operator Pick class, but to the smaller {\em operator Loewner class}.
\begin{proposition}\label{Mloewner}
With the notation of Definition {\rm \ref{defMatResolv}}, there exists an analytic operator-valued function
$F:\Pi^n \to \mathcal{L}(\h)$ such that
for all $z,\ w\in\Pi^n$,
\beq\label{loewDecomp}
M(z)-M(w)^* = F(w)^*(z-\bar w)_PF(z)
\eeq
on $\h$.
\end{proposition}
\begin{proof}
The identity \eqref{imM} shows that such a relation holds on $\N\oplus \mathcal{D}(A)$; we must extend it to all of $\h$.
Write $P_j$ as an operator matrix with respect to the decomposition $\h=\N\oplus\M$, as in equation \eqref{expP}.  Then $z_P$ has the matricial expression \eqref{formzP}.  For $z\in\Pi^n$ let
\[
F^\sharp(z) = \left(\bbm 1&0\\0&A\ebm - \bbm 0&0\\0&1\ebm z_P\right)\inv \bbm i&0\\0&1+iA\ebm.
\]
Then $F^\sharp(z)$ is an operator from $\N\oplus\mathcal{D}(A)$ to $\h$, and we find that
\begin{align*}
F^\sharp(z) &= \bbm 1& 0\\-z_{B^*} & A-z_Y\ebm\inv \bbm i&0\\0&1+iA\ebm \\
	&=\bbm i&0 \\ i(A-z_Y)\inv z_{B^*} & (A-z_Y)\inv (1+iA)\ebm : \N\oplus\mathcal{D}(A) \to \h.
\end{align*}
Let
\beq\label{defFz}
F(z) = \bbm i&0 \\ i(A-z_Y)\inv z_{B^*} & i+(A-z_Y)\inv (1+iz_Y)\ebm : \N\oplus\M \to \h.
\eeq
Since
\[
(A-z_Y)\inv (1+iA)=i+ (A-z_Y)\inv (1+iz_Y)
\]
on $\N\oplus\mathcal{D}(A)$ and the right hand side of the last equation is a bounded operator on all of $\h$, it is clear that, for every $z\in\Pi^n$, $F(z)$ is a continuous extension to $\h$ of $F^\sharp(z)$ and is a bounded operator.   Furthermore $F$ is analytic on $\Pi^n$.

By Proposition \ref{alternM}, equation \eqref{imM}, the relation \eqref{loewDecomp} holds on the dense subspace $\N\oplus\mathcal{D}(A)$ of $\h$ for every $z,\ w\in\Pi^n$.  Since the operators on both sides of equation \eqref{loewDecomp} are continuous on $\h$, the equation holds throughout $\h$.
\end{proof}
\begin{corollary}\label{MatPick}
A matricial resolvent has a non-negative imaginary part at every point of $\Pi^n$.
\end{corollary}
\begin{proof}
In the notation of Proposition \ref{Mloewner}, on choosing $w=z$ in equation \eqref{loewDecomp} and dividing by $2i$ we obtain the relation
\[
\im M(z)= F(z)^*(\im z_P )F(z)
\]
on $\h$.    We have
\[
\im z_P =\sum_j (\im z_j)P_j  \geq 0,
\]
and so $\im M(z) \geq 0$ on $\h$ for all $z\in\Pi^n$.
\end{proof}

Here is a concrete example of a matricial resolvent.
\begin{example}\label{type4MatR}
\rm
The function
\beq\label{formtype4}
M(z)= \frac{1}{z_1+z_2}\bbm 2z_1z_2 & i(z_1-z_2) \\ -i(z_1-z_2) & -2 \ebm
\eeq
is the matricial resolvent corresponding to
\[
\h=\C^2, \quad \N=\M=\C,\quad A=0 \mbox{ on }\C,\quad P_1= \half \bbm 1&1\\1&1 \ebm, \quad P_2=1-P_1.
\]
\end{example}

\section{Nevanlinna representations of type $4$} \label{type4}
In this section we derive a multivariable analog of the most general form of Nevanlinna representation for functions in the one-variable Pick class (Theorem \ref{thm1.1}).  We start with a multivariable Herglotz theorem \cite[Theorem 1.8]{ag90}.  We shall say that an analytic operator-valued function $F$ on $\D^n$ is a {\em  Herglotz function} if $\re F(\la) \geq 0$ for all $\la\in\D^n$.  For present purposes we need the following modification of the notion.
\begin{definition}\label{strongHerg}
An analytic function  $F: \D^n \to \mathcal {L(K)}$, where $\mathcal K$ is a Hilbert space, is a  {\em strong Herglotz function} if, for every commuting $n$-tuple $T=(T_1,\dots,T_n)$ of operators on a Hilbert space and for $0\leq r < 1, \,  \re F(rT) \geq 0$.  
\end{definition}
 In \cite{ag90} these functions were called $\mathcal{F}_n$-Herglotz functions.  The  class of strong Herglotz functions has also been called the {\em Herglotz-Agler class} (for example \cite{kaluzh,BKV2}).
It is clear that every strong Herglotz function is a Herglotz function, and in the cases $n=1$ and $2$ the converse is also true \cite{ag90}.

\begin{theorem}\label{thm2.1}
Let $\mathcal K$ be a Hilbert space and let $F: \D^2 \to \mathcal {L(K)}$ be a strong Herglotz function such that $F(0)=1$.  There exist a Hilbert space $\mathcal H$, an orthogonal decomposition $P$ of $\mathcal H$, an isometric linear operator $V: \mathcal K \to \mathcal H$ and a unitary operator $U$ on $\mathcal H$ such that, for all $\la\in\D^n$,
\beq\label{aglerRep}
F(\lambda) = V\ad \frac{1+U\lambda_P}{1 - U\lambda_P} V.
\eeq

Conversely, every function $F: \D^n \to \mathcal{ L( K)}$ expressible in the form \eqref{aglerRep} for some $\mathcal{H},\ P,  \ V$ and $U$ with the stated properties  is a strong Herglotz function and satisfies $F(0)=1$.
\end{theorem}

Note that  $\la_P=\sum_j \la_j P_j$ has operator norm at most $\|\la\|_\infty < 1$  for $\la\in\D^n$, and hence equation \eqref{aglerRep} does define $F$ as an analytic operator-valued function on $\D^n$.

On specialising to scalar-valued functions in the $n$-variable Herglotz class we obtain the following consequence.
\begin{corollary}\label{thm2.2}
Let $f$ be a scalar-valued strong Herglotz function on $\D^n$. There exists a Hilbert space $\mathcal H$, a unitary operator $L$ on $\mathcal H$, an orthogonal decomposition $P$ of $\mathcal H$,  a real number $a$ and a vector $v \in \mathcal H$   such that, for all $\la\in\D^n$,
\beq\label{scalarHerg}
 f(\lambda) = -ia + \ip{(L-\lambda_P)^{-1}(L+\lambda_P)v}{v}.
\eeq
 Conversely, for any $ \mathcal{H}, L, P, a$  and $v$ with the properties described, equation \eqref{scalarHerg} defines $f$ as an $n$-variable strong Herglotz function.  
\end{corollary}

Again, the right hand side of equation \eqref{scalarHerg} is an analytic function of $\la\in\D^n$ since
\[
(L-\lambda_P)^{-1}= L^{-1}(1-\la_P L^{-1})^{-1}
\]
is a bounded operator and is analytic in $\la$.

\begin{definition}\label{def2.1}
A {\em Nevanlinna representation of type $4$} of a function $h:\Pi^n\to\C$ consists of an orthogonally decomposed Hilbert space $\mathcal {H=N\oplus M}$,  a self-adjoint densely defined operator $A$ on $\mathcal {M}$, an orthogonal decomposition $P$ of $\h$, a real number $a$ and a vector $v \in \mathcal H$ such that 
\beq\label{formh}
h(z)=a+\ip{ M(z)v}{v}
\eeq
for all $z\in\Pi^n$, where $M(z)$ is the structured resolvent of $A$ of type $4$ corresponding to $P$  (given by the formula \eqref{defM}).
 \end{definition}

We wish to convert Corollary \ref{thm2.2} to a representation theorem for suitable analytic functions on $\Pi^n$.  The fact that the corollary only applies to {\em strong} Herglotz functions results in representation theorems for a subclass of the Pick class $\Pick_n$.  Recall from the introduction:
\begin{definition}\label{loewner}
The {\em Loewner class} $\p$ is the set of analytic functions $h$ on $\Pi^n$ with the property that there exist $n$ positive semi-definite functions $A_1,\dots,A_n$ on $\Pi^n$, analytic in the first argument,  such that
\[
h(z) - \overline{h(w)} = \sum_{j=1}^n (z_j-\overline{w_j}) A_j(z,w)
\]
for all $z,w\in\Pi^n$.
\end{definition}
A function $h$ on $\Pi^n$ belongs to $\p$ if and only if it corresponds under conjugation by the Cayley transform to a function in the Schur-Agler class of the polydisc \cite[Lemma 2.13]{AMY}.  Another characterization: $h\in\p$ if and only if, for every commuting $n$-tuple $T$ of bounded operators with strictly positive imaginary parts, $h(T)$ has positive imaginary part.

We can now prove Theorem \ref{thm2.4} from the introduction:  {\em a function $h$  defined on $\Pi^n$ has a Nevanlinna representation of type $4$ if and only if $h \in \p$.}
\begin{proof}
Let $h\in\pick$.  Define an $n$-variable Herglotz function $f:\D^n \to\C$ by
\beq\label{deff}
f(\la) =-i h(z)
\eeq
where
\beq
z_j= i\frac{1+\la_j}{1-\la_j} \qquad \mbox{ for } j=1,\dots,n.
\eeq
When $\la\in \D^n$ the point $z$ belongs to $\Pi^n$, and so $f(\la) $ is well defined, and since $\im h(z) \geq 0$ we have $\re f(\la) \geq 0$, so that $f$ is indeed a Herglotz function.  In fact $f$ is even a strong Herglotz function: since $h\in\p$, the function $\ph\in\s$ corresponding to $h$ lies in the Schur-Agler class of the polydisc, and so $f=(1+\ph)/(1-\ph)$ is a strong Herglotz function.

By Corollary \ref{thm2.2} there exist a real number $a$, a Hilbert space $\mathcal H$,  a vector $v \in \mathcal H$, a unitary operator $L$ on $\h$ and an orthogonal decomposition $P$ on $\h$ such that, for all $z\in\Pi^n$,
\begin{align}\label{hif}
 h(z) &= if(\la)= a+ \ip{i(L-\la)^{-1}(L+\la)v}{v} \nn \\
	& = a + \ip{i[L - (z-i)(z+i)^{-1}]^{-1}[L + (z-i)(z+i)^{-1}]v}{v}.
\end{align}
Here and in the rest of this section $z,\ \la$ are identified with the operators $z_P,\ \la_P$ on $\h$, and in consequence the relation
\[
\la= \frac{z-i}{z+i}
\]
is meaningful and valid.

For $z\in\Pi^n$ let
\beq\label{defML}
 M(z)=  i\left(L - \la\right)^{-1}\left(L + \la\right)= i\left(L - \frac{z-i}{z+i}\right)^{-1}\left(L + \frac{z-i}{z+i}\right).
\eeq
Since $L$ is unitary on $\h$ and $\la\in\D^n$, the operator $M(z)$ is bounded on $\h$ for every $z\in\Pi^n$
and, by equation \eqref{hif}, we have
\beq\label{repofh}
h(z)=a+\ip{ M(z)v}{v}
\eeq
for all $z\in\Pi^2$.    Theorem \ref{thm2.4} will follow provided we can show that  $ M(z)$ is given by equation \eqref{defM} for a suitable self-adjoint operator $A$.

Observe that
\begin{align}\label{exp2.1}
  M(z) &= i((z+i)L - (z-i))\inv ((z+i)L + (z-i)) \notag \\
 &= i\left(z(L-1) + i(L+1)\right)\inv\left(z(L+1) +i (L-1)\right).
\end{align}
We wish to take out a factor $1-L$ from both factors in equation \eqref{exp2.1}, but this may be impossible since $1-L$ can have a nonzero kernel.  Accordingly we decompose $\h$ into $\mathcal{N}\oplus\mathcal{M}$ where $\mathcal{N} = \ker (1-L), \  \M= \N^\perp$.  With respect to this decomposition we can write 
 $L$  as an operator matrix
\[
 L = \begin{bmatrix}
      1 & 0 \\ 0 & L_0
     \end{bmatrix},
\]
where $L_0$ is unitary and $\ker (1 - L_0) = \{0\}$. Substituting into equation \eqref{exp2.1} we have
\begin{align} \label{exp2.1.1}
 M(z)  &= i\left(z\begin{bmatrix}
      0 & 0 \\ 0 & L_0-1
     \end{bmatrix} + 
  i\begin{bmatrix}
  2 & 0 \\ 0 & L_0 +1
 \end{bmatrix}\right)\inv \left(
  z\begin{bmatrix}
  2 & 0 \\ 0 &  L_0 +1
 \end{bmatrix} + 
  i\begin{bmatrix}
  0 & 0 \\ 0 & L_0 -1
 \end{bmatrix} z\right)  \nn \\
	&=\left(-z \bbm 0&0 \\0& 1-L_0 \ebm +\bbm 2i & 0 \\ 0 & i(1+L_0)\ebm \right)\inv
	 \left( z\bbm 2i &0 \\ 0&i(1+L_0) \ebm + \bbm 0&0\\0 & 1-L_0\ebm \right)
\end{align}
 Formally we may now write
\begin{align} \label{newtilde}
 M(z) &= \begin{bmatrix}
     -\half i & 0 \\ 0 & (1 - L_0)\inv
    \end{bmatrix}
\left( -z \begin{bmatrix}
      0 & 0 \\ 0 &1
     \end{bmatrix}  + 
 \begin{bmatrix}
  1 & 0 \\ 0 &  i\frac{1 + L_0}{1 - L_0}
 \end{bmatrix}\right)\inv  \times  \nn \\
	& \hspace{2cm}  \left(
  z\begin{bmatrix}
  1 & 0 \\ 0 & i\frac{1 + L_0}{1 - L_0}
 \end{bmatrix} +
 \begin{bmatrix}
  0 & 0 \\ 0 & 1
 \end{bmatrix} \right)
 \begin{bmatrix}
  2 i & 0 \\ 0 & 1 - L_0
 \end{bmatrix},
\end{align}
but whereas equation \eqref{exp2.1.1} is a relation between bounded operators defined on all of $\h$, equation \eqref{newtilde} involves unbounded, partially defined operators and we must verify that the product of operators on the right hand side is meaningful.

Let
\[
 A = i\frac{1 + L_0}{1 - L_0}.
\]
Since $L_0$ is unitary on $\mathcal{M}$ and $\ker (1-L_0) =\{0\}$, the operator $A$ is self-adjoint and densely defined on $\mathcal{M}$ \cite[Section 121]{RN}.   The domain $\mathcal{D}(A)$ of $A$ is the dense subspace $\ran (1-L_0)$ of $\M$.   It follows from the definition of $A$ that
\beq\label{useful}
 (1-L_0)^{-1} =  \half (1-iA),
\eeq
which is an equation between bijective operators from $\mathcal{D}(A)$ to $\M$.  Likewise
\beq\label{useful2}
1+L_0 = -2iA(1-iA)^{-1} :\M\to\mathcal{D}(A)
\eeq
are bounded operators.

Let us continue the calculation from the first factor on the right hand side of equation \eqref{exp2.1.1}.   Since $\ker (1-L_0) = \{0\}$, the right hand side of the relation
\begin{align*}
-z \bbm 0&0\\0& 1-L_0  \ebm+ \bbm 2i&0\\0& i(1+L_0)\ebm &= \left( -z \bbm 0&0\\0&1 \ebm + \bbm 1 &0\\ 0&A \ebm \right) \bbm 2i&0\\0& 1-L_0 \ebm  
\end{align*}
comprises a bijective map from $\h$ to $\N \oplus \mathcal D(A)$ followed by a bijection from $\N \oplus \mathcal D(A)$ to $\h$ (recall the equation \eqref{needThis}).
We may therefore take inverses in the equation to obtain
\begin{align}\label{2nd}
\left(    -z \bbm 0&0\\0& 1-L_0  \ebm+ \bbm 2i&0\\0& i(1+L_0)\ebm \right)^{-1} 
	&= \bbm -\half i&0\\0&(1-L_0)^{-1} \ebm \left( \bbm1&0\\0&A \ebm - z\bbm 0&0\\0&1 \ebm \right)^{-1}  \nn \\
	&=  \bbm -\half i&0\\0& \half (1-iA) \ebm \left( \bbm1&0\\0&A \ebm - z\bbm 0&0\\0&1 \ebm \right)^{-1} 
\end{align}
as operators on $\N\oplus \mathcal{D}(A)$.	

Similar reasoning applies to the equation
\begin{align}\label{1st}
z\bbm 2i & 0\\0& i(1+L_0) \ebm  + \bbm 0&0\\0& 1-L_0 \ebm  &= \left(z\bbm 1&0\\0& A \ebm + \bbm0&0\\0&1 \ebm \right) \bbm 2i&0\\ 0&1-L_0 \ebm  \nn\\
   &=\left(z\bbm 1&0\\0& A \ebm + \bbm0&0\\0&1 \ebm \right)  \bbm -\half i &0\\0&\half (1-iA) \ebm\inv;  
\end{align}
it is valid as an equation between operators on $\h$.  The right hand side comprises an operator from $\h$ to $\N\oplus \mathcal{D}(A)$ followed by an operator from $\N\oplus \mathcal{D}(A)$ to $\h$,
 and so both sides of the equation denote an operator on $\h$.  

On combining equations \eqref{exp2.1.1}, \eqref{2nd} and \eqref{1st} we obtain
\begin{align*}
M(z) &=  \bbm -\half i&0\\0& \half (1-iA) \ebm \left( \bbm1&0\\0&A \ebm - z\bbm 0&0\\0&1 \ebm \right)^{-1} \left(z\bbm 1&0\\0& A \ebm + \bbm0&0\\0&1 \ebm \right)  \bbm -\half i &0\\0&\half (1-iA) \ebm\inv.
\end{align*}
Premultiply this equation by $2$ and postmultiply by $\half$ to deduce that $M(z)$ is indeed the structured resolvent of $A$ of type $4$ corresponding to $P$, as defined in equation \eqref{defM}.  Thus the formula \eqref{repofh} is a Nevanlinna representation of $h$ of type $4$.

Conversely,  let $h\in\p$ have a type 4 representation \eqref{formh}. 
By Proposition \ref{Mloewner} there exists an analytic operator-valued function
$F:\Pi^n \to \mathcal{L}(\h)$ such that,
for all $z,\ w\in\Pi^n$,
\beq
M(z)-M(w)^* = F(w)^*(z-\bar w)_PF(z)
\eeq
on $\h$.  Hence
\begin{align*}
h(z)-\overline{h(w)} &= \ip{ (M(z)-M(w)^*)v}{v} \\
	&=\ip{F(w)^*(z-\bar w)_P F(z)v}{v} \\
	&=\sum_{j=1}^n (z_j- \bar w_j)A_j(z,w)
\end{align*}
for all $z, \ w\in\Pi^n$, where
\[
A_j(z,w) = \ip{P_jF(z)v}{F(w)v}.
\]
The $A_j$ are clearly positive semidefinite on $\Pi^n$, and hence $h$ belongs to the Loewner class $\mathcal{L}_n$.
\end{proof}

\section{Nevanlinna representations of types 3, 2 and 1}\label{type321}
Nevanlinna representations of type $4$ have the virtue of being general for functions in $\p$, but they are undeniably cumbersome.
In this section we shall show that there are three simpler representation formulae, corresponding to increasingly stringent growth conditions on $h \in \pick$.

In Nevanlinna's one-variable representation formula of Theorem \ref{thm1.1},
  \beq\label{OneVar}
   h(z) = a+bz+\int \frac{1+tz}{t-z} \  \dd\mu(t),
  \eeq
it may be the case for a particular $h\in\Pick$ that the $bz$ term is absent.  The analogous situation in two variables is that the space $\N$ in a type 4 representation may be zero.  Equivalently, in the corresponding Herglotz representation, the unitary operator $L$ does not have $1$ as an eigenvalue.  This suggests the following notion. 
\begin{definition}\label{def3.1}
A  {\em Nevanlinna representation of type} $3$ of a function $h$ on $\Pi^n$ consists of a Hilbert space $\mathcal H$, a self-adjoint densely defined operator $A$ on $\h$, a positive decomposition $Y$ of $\h$, a real number $a$ and a vector $v \in \mathcal H$ such that, for all $z\in\Pi^n$,
\beq\label{type3rep}
h(z) = a+ \ip{(1-iA)(A - z_Y)^{-1} (1+z_YA)(1-iA)\inv v}{v}.
\eeq
\end{definition}
Thus $h$ has a type 3 representation if  $h(z) = a+ \ip{ M(z)v}{v}$ where $ M(z)$ is the structured resolvent of $A$ of type $3$ corresponding to $Y$, as given by equation \eqref{defM3}.

In \cite{ATY} the authors derived a somewhat simpler representation which can also be regarded as an analog of the case $b=0$ of Nevanlinna's one-variable formula \eqref{OneVar}. 
\begin{definition} \label{defType2ref}
A {\em Nevanlinna representation of type} $2$ of a function $h$ on $\Pi^n$ consists of a Hilbert space $\mathcal H$, a self-adjoint densely defined operator $A$ on $\h$, a positive decomposition $Y$ of $\h$, a real number $a$ and a vector $\alpha \in \mathcal H$ such that, for all $z\in\Pi^n$
\beq \label{type2rep}
h(z)=a+\ip{(A-z_Y)\inv \al}{\al}.
 \eeq
\end{definition}
This means of course that, for all $z\in\Pi^n$,
\[
h(z) = a + \ip{M(z)\alpha}{\alpha}
\]
where $M(z)$ is the structured resolvent of $A$ of type $2$ corresponding to $Y$ (compare equation \eqref{resolv2}).

We wish to understand the relationship between type 3 and type 2 representations.
\begin{proposition}\label{2gives3}
 If $h \in \mathcal P_n$  has a type $2$ representation then $h$ has a type $3$ representation.   Conversely,
 if $h \in \mathcal P_n$ has a type $3$ representation as in equation \eqref{type3rep} with the additional property that $v \in \mathcal D(A)$ then $h$ has a type $2$ representation. 
\end{proposition}
\begin{proof}
 Suppose that $h \in \Pick_n$ has the type 2 representation 
\[
 h(z) = a_0 + \ip{(A - z_Y)^{-1}\alpha}{\alpha}
\]
for some $a_0\in\R$, positive decomposition $Y$ and $\al\in \h$.
We must show that $h$ has a representation of the form \eqref{type3rep}
for some $a\in\R$ and $v\in \h$.   By Proposition \ref{2ndResolv3}, it suffices to find $a\in\R$ and $v\in\mathcal D(A)$ such that
\[
h(z) = a+ \ip{(1-iA)\left\{(A - z_Y)^{-1}-A(1+A^2)\inv\right\}(1+iA)v}{v}
\]
for all $z\in\Pi^n$.

To this end, let $C=A(1+A^2)\inv$ and let
\beq\label{aanda_0}
a = a_0 + \ip{C\alpha}{\alpha}.
\eeq
 Since  $1+iA$ is invertible on $\h$ and  $\ran (1+iA)^{-1}\subset \mathcal D(A)$ we may define
\beq\label{vandalpha}
 v = (1+iA)^{-1} \alpha \in \mathcal D(A).
\eeq

Then
\begin{align*}
 h(z) &= a_0 + \ip{(A - z_Y)^{-1}\alpha}{\alpha} \notag \\
	&= a-\ip{C\al}{\al} + \ip{(A - z_Y)^{-1}\alpha}{\alpha} \notag \\
 	&= a+ \ip{\left\{(A - z_Y)^{-1}-C\right\}(1+iA)v}{(1+iA)v} \notag \\
 	&= a+ \ip{(1-iA)\left\{(A - z_Y)^{-1}-C\right\}(1+iA)v}{v}
\end{align*}
as required.  Thus $h$ has a type $3$ representation.

Conversely, let $h$ have a type 3 representation \eqref{type3rep} such that $v \in \mathcal D(A)$, that is
\[
h(z)=a+\ip{M(z)v}{v}
\] 
where $a\in\R$ and $M$ is the structured resolvent of $A$ of type 3 corresponding to $Y$, as in equation \eqref{defM3}.
  Since $v \in \mathcal D(A)$ we may define the vector $\alpha \df (1+iA)v \in \h$, and furthermore, by Proposition \ref{2ndResolv3},
\begin{align*}
 h(z) &= a + \ip{(1-iA)\left\{(A - z_Y)^{-1} - C\right\}(1+iA)v}{v} \\
 	&= a + \ip{\left\{(A - z_Y)^{-1} - C\right\}\alpha}{\alpha}  \\
 	&= a - \ip{C\alpha}{\alpha} + \ip{(A - z_Y)^{-1}\alpha}{\alpha}  \\
 	&= a_0 + \ip{(A - z_Y)^{-1}\alpha}{\alpha} ,
\end{align*}
where $a_0\in\R$ is given by equation \eqref{aanda_0}. Thus $h$ has a representation of type 2.
\end{proof}

A special case of a type 2 representation occurs when the constant term $a$ in equation \eqref{type2rep} is $0$. In one variable, this corresponds to Nevanlinna's  characterization of the Cauchy transforms of positive finite measures on $\R$.
Accordingly we define a  {\em type $1$ representation} of $h\in \p$ to be the special case of a type 2 representation of $h$ in which $a=0$ in \eqref{type2rep}.  

\begin{definition}\label{def4.1}
 An analytic function $h$ on $\Pi^n$ has a {\em Nevanlinna representation of type $1$}  if there exist a Hilbert space $\h$, a densely defined self-adjoint operator $A$ on $\h$, a positive decomposition $Y$ of $\h$  and  a vector $\alpha\in\h$ such that, for all $z\in\Pi^n$,
 \beq\label{type1formula}
  h(z) = \ip{(A - z_Y)^{-1}\alpha}{\alpha}.
 \eeq
\end{definition}

A representation of type 1 is obviously a representation of type 2. The following proposition is an immediate corollary of Proposition \ref{2gives3}.

\begin{proposition}\label{1equiv3}
 A function $h \in \pick$  has a type $1$ representation if and only if $h$ has a type $3$ representation as in equation \eqref{type3rep} with the additional properties that $v \in \mathcal D(A)$ and
\[
 a - \ip{A(1+A^2)\inv\alpha}{\alpha} =0.
\]
\end{proposition}

For consistency with our earlier terminology for structured resolvents and representations we should have to define a structured resolvent of type 1 to be the same as a structured resolvent of type 2.  We refrain from making such a confusing definition.

We conclude this section by giving examples of the four types of Nevanlinna representation in two variables.
\begin{example}\label{4types}
\rm
(1) The formula
\[
h(z)= -\frac{1}{z_1+z_2} = \ip{(0-z_Y)^{-1}v}{v}_{\C},
\]
where $Y=(\half,\half)$ and $v=1/\sqrt{2}$, exhibits a representation of type $1$, with $A=0$.\\

\nin (2)  Likewise
\[
h(z)= 1-\frac{1}{z_1+z_2} =1+ \ip{(0-z_Y)^{-1}v}{v}_{\C}
\]
is a representation of type $2$.\\

\nin (3)  Let
\beq\label{trueh}
h(z) =\threepartdef{ \ds  \frac{1}{1+z_1z_2}\left(z_1-z_2 + \frac{iz_2(1+z_1^2)}{\sqrt{z_1z_2}}\right)}{z_1z_2 \neq -1}{}{}{\half (z_1+z_2)}{z_1z_2 = -1}
\eeq
where we take the branch of the square root that is analytic in $\C\setminus [0,\infty)$ with range $\Pi$.
We claim that $h\in\Pick_2$ and that $h$ has the type $3$ representation
\beq\label{temp}
h(z)= \ip{M(z)v}{v}_{L^2(\R)},
\eeq
where $M(z)$ is the structured resolvent of type $3$ given in Example \ref{type3ex2} and $v(t)=1/\sqrt{\pi(1+t^2)}$.     To see this, let $h$ be temporarily defined by equation \eqref{temp}.  Since $v$ is an even function in $L^2(\R)$, equation \eqref{fcnoftype3} tells us that 
\[
h(z) = \int_{-\infty}^\infty \frac {t(1+z_1z_2)+(1-it)(itz_1+z_2)}{\pi(t^2-z_1z_2)(1+t^2)} \ \dd t.
\]
Since the denominator is an even function of $t$, the integrals of all the odd powers of $t$ in the numerator vanish, and we have, provided $z_1z_2\neq -1$,
\begin{align*}
h(z) &= \frac{2}{\pi}\int_{0}^\infty \frac {z_2 +t^2z_1}{(t^2-z_1z_2)(1+t^2)} \ \dd t\\
	&=  \frac{2}{\pi}\int_{0}^\infty \frac {z_2(1+z_1^2)}{1+z_1z_2}\, \, \frac{1}{t^2 -z_1z_2} + \frac{z_1-z_2}{1+z_1z_2} \, \,\frac{1}{1+t^2} \ \dd t.
\end{align*}
Now, for $w\in\Pi$,
\[
\int_0^\infty \frac{\dd t}{t^2-w^2} = \frac{i\pi}{2w},
\]
and so we find that $h$ is indeed given by equation \eqref{trueh} in the case that $z_1z_2 \neq -1$.  When $z_1z_2=-1$ we have
\begin{align*}
h(z) &= \frac{2}{\pi}\int_0^\infty \frac{z_2+z_1t^2}{(1+t^2)^2} \dd t \\
	&= \frac{2}{\pi}\int_0^\infty \frac{z_1}{1+t^2} + \frac{z_2-z_1}{(1+t^2)^2} \dd t\\
	&= \half(z_1+z_2).
\end{align*}
 Thus equation \eqref{temp} is a type $3$ representation of the function $h$ given by equation \eqref{trueh}.  This function is {\em constant} and equal to $i$ on the diagonal $z_1=z_2$.\\

\nin (4)  The function
\[
h(z)= \frac{z_1 z_2}{z_1+z_2} = -\left(-\frac{1}{z_1}-\frac{1}{z_2}\right)^{-1}
\]
clearly belongs to $\Pick_2$.  It has the representation of type $4$
\[
h(z)= \ip{M(z)v}{v}_{\C^2}
\]
where $M(z)$ is the matricial resolvent given in Example \ref{type4MatR} and
\[
 v=\frac{1}{\sqrt{2}}\begin{pmatrix} 1 \\ 0 \end{pmatrix}.
\]
\end{example}
We claim that each of the above representations is of the simplest available type for the function in question; for example, the function $h$ in part (4) does not have a Nevanlinna representation of type 3.  To prove this claim (which we shall do in Example \ref{fourtypes} below) we need characterizations of the types of functions  -- the subject of the next two sections.

\section{Asymptotic behavior and types of representations}\label{asymptotic}

In this section we shall give function-theoretic conditions for a function  in $ \pick$ to have a representation of a given type. These conditions will be in terms of the asymptotic behavior of the function at $\infty$. 

Every function in $\pick$ has a type $4$ representation, by Theorem \ref{thm2.4}.  Let us characterize the functions that possess a type $3$ representation.  We denote by $\chi$ the vector $(1,\dots,1)$ of ones in $\C^n$.
The following statement contains Theorem \ref{type3intro}.
\begin{theorem}\label{type3asymp}
 The following three conditions are equivalent for a function $h \in \pick$.
\begin{enumerate}
\item  The function $h$ has a Nevanlinna representation of type $3$;
\item \beq\label{type3liminf}
\liminf_{s\to\infty} \frac{1}{s}  \im h(is\chi)=0;
\eeq
 \item \beq\label{type3lim}
\lim_{s\to\infty} \frac{1}{s} \im h(is\chi)=0.
\eeq
\end{enumerate}
\end{theorem}
\begin{proof}
 (1)$\Rightarrow$(3)  Suppose that $h$ has a Nevanlinna representation of type $3$:
\beq\label{type3bis}
h(z)= a + \ip{(1-iA)(A-z_Y)\inv(1+z_YA)(1-iA)\inv v}{v}
\eeq
for suitable $a\in\R, \h, A, Y$ and $v\in\h$.
 Since
\[
 (is\chi)_Y =\sum_j isY_j  = is
\]
we have
\[
 h(is\chi) =a+\ip{(1-iA)(A-is)\inv(1+isA)(1-iA)\inv v}{v}.
\]

Let $\nu$ be the scalar spectral measure for $A$ corresponding to the vector $v\in\h$.  By the Spectral Theorem
\begin{align*}
h(is\chi)&=a+\int (1-it)(t-is)\inv(1+ist)(1-it)\inv \ \dd\nu(t)\\
	&= a+\int \frac {1+ist}{t-is} \ \dd\nu(t).
\end{align*}
Since
\[
\im\frac{1+ist}{t-is} = \frac{s(1+t^2)}{s^2+t^2},
\]
we have
\[
\frac{1}{s} \im h(is\chi) = \int \frac{1+t^2}{s^2+t^2} \ \dd\nu(t).
\]
The integrand decreases monotonically to $0$ as $s\to\infty$ and so, by the Monotone Convergence Theorem,  equation \eqref{type3lim} holds.

(3)$\Rightarrow$(2) is trivial.

(2)$\Rightarrow$(1)
Now suppose that $h \in \pick$ and 
\[
\liminf_{s\to\infty}\frac{1}{s}\IM h(is\chi)=0.
\]
 By Theorem \ref{thm2.4}, $h$ has a Nevanlinna representation of type $4$: that is, there exist $a, \h, \mathcal{N}\subset\h$,  operators $A, \, Y$ on $\mathcal{N}^\perp$ and a vector $v\in\h$ with the properties described in Definition \ref{def3.1} such that
\[
h(z) = a+ \ip{M(z)v}{v}
\]
for all $z\in\Pi^n$, where
\beq\label{defM4bis}
M(z)=\bbm -i&0\\0&1-iA \ebm
   \left( \bbm 1&0 \\ 0 & A \ebm  -
   z_P\bbm 0&0\\0& 1 \ebm \right)\inv 
\left(z_P\bbm 1&0\\0& A\ebm +
   \bbm 0&0\\ 0&1 \ebm \right)
   \bbm  -i&0\\0&1-iA \ebm^{-1}.
\eeq
Thus, for $s>0$, since once again $(is\chi)_P=is$,
\begin{align*}
 M(is\chi) &= \bbm -i&0\\0&1-iA \ebm
    \bbm 1&0 \\ 0 & (A-is)\inv \ebm 
   \bbm is&0\\0& 1+isA \ebm 
   \bbm  i&0\\0&(1-iA)\inv \ebm\\
	&=\bbm is&0\\0& (1-iA)(A-is)\inv(1+isA)(1-iA)\inv\ebm.
\end{align*}

Let the projections of $v$ onto $\mathcal{N},\ \mathcal{N}^\perp$ be $v_1, v_2$ respectively.  Then 
\begin{align*} 
h(is\chi) &=  a+\ip{M(is\chi)v}{v}  \\
	&=  a  +is\norm{v_1}^2 + \ip{(1-iA)(A-is)\inv(1+isA)(1-iA)\inv v_2}{v_2} 
\end{align*}
and therefore
\begin{align*}
\frac{1}{s} \im h(is\chi) &=  \norm{v_1}^2 + \frac{1}{s} \im\ip{(1-iA)(A-is)\inv(1+isA)(1-iA)\inv v_2}{v_2}\\
	&\geq \norm{v_1}^2 
\end{align*}
by Corollary \ref{Mpick}.  Hence
\begin{align*}
0 &= \liminf_{s\to\infty}\frac{1}{s}\im h(is\chi) \\
	& \geq \norm{v_1}^2. 
\end{align*}
It follows that $v_1 = 0$.

  Let the compression of the projection $P_j$ to $\mathcal{N}^\perp$ be $Y_j$: then $Y=(Y_1,\dots,Y_n)$ is a positive decomposition of $\mathcal{N}^\perp$, and  the compression of $z_P$ to $\mathcal{N}^\perp$ is  $z_Y$.  
By Remark \ref{3.2} the (2,2) block $M_{22}(z)$ in $M(z)$ is
\begin{align*}
M_{22}(z)&= (1-iA)(A-z_Y)\inv (1+z_YA)(1-iA)\inv.
\end{align*}
Since $v_1=0$ it follows that
\begin{align*}
h(z)&= a+ \ip{M(z)v}{v}   \\
	&= a+ \ip{M_{22}(z)v_2}{v_2} \\
	&= a+ \ip{(1-iA)(A-z_Y)\inv (1+z_YA)(1-iA)\inv v_2}{v_2},
\end{align*}
which is the desired type $3$ representation of $h$.   Hence (2)$\Rightarrow$(1).
\end{proof}
In \cite{BKV2} it is shown that condition (3) in the above theorem is also a necessary and sufficient condition that $-ih$ have a {\em $\Pi^n$-impedance-conservative realization.}

Type $2$ representations were characterized by the following theorem in \cite{ATY} in the case of two variables.  The following result, which contains Theorem \ref{type2intro}, shows that the result holds generally.
\begin{theorem}\label{type2asymp}
 The following three conditions are equivalent for a function $h\in\pick$.
\begin{enumerate}
\item The function $h$ has a Nevanlinna representation of type $2$;
\item \beq \label{type2liminf}\liminf_{s\to\infty} s\IM h(is\chi)<\infty;\eeq
\item \beq \label{type2lim} \lim_{s\to\infty} s\IM h(is\chi)<\infty.\eeq
\end{enumerate}
\end{theorem}
\begin{proof}
\nin (1)$\Rightarrow$(3)  Suppose that $h$ has the type $2$ representation $h(z)=a+\ip{(A-z_Y)\inv v}{v}$ for a suitable real $a$, self-adjoint $A$, positive decomposition $Y$ and vector $v$.  Let $\nu$ be the scalar spectral measure for $A$ corresponding to the vector $v$.  Then, for $s>0$, $A-(is\chi)_Y = A-is$ and so
\begin{align*}
s \im h(is\chi) &= s \im \int \frac{\dd \nu(t)}{t-is} \\
	&= \int \frac{ s^2 \ \dd\nu(t)}{t^2+s^2}.
\end{align*}
The integrand is positive and increases monotonically to $1$ as $s\to\infty$.  Hence, by the Dominated Convergence Theorem
\[
\lim_{s\to\infty} s \im h(is\chi) = \nu(\R) = \|v\|^2 < \infty.
\]
Hence (1)$\Rightarrow$(3).

\nin (3)$\Rightarrow$(2) is trivial.

\nin (2)$\Rightarrow$(1)  Suppose (2) holds.  {\em A fortiori},
\[
\liminf_{s\to\infty} \frac{1}{s} \im h(is\chi) = 0.
\]
By Theorem \ref{type3asymp} $h$ has a type $3$ representation \eqref{type3bis}
for suitable $a\in\R, \h, A, Y$ and $v\in\h$. Let $\nu$ be the scalar spectral measure for $A$ corresponding to the vector $v$.  Then for $s>0$
\begin{align*}
s\im h(is\chi) &= s\im \int\frac{1+ist}{t-is} \ \dd\nu(t) \\
	&= \int \frac{s^2(1+t^2)}{t^2+s^2}\ \dd\nu(t). 
\end{align*}
As $s\to\infty$ the integrand increases monotonically to $1+t^2$.  Condition (2) now implies that
\[
\int 1+t^2 \ \dd\nu(t) < \infty.
\]
It follows that $v\in\mathcal{D}(A)$.  Hence, by Proposition \ref{2gives3}, $h$ has a representation of type $2$.
\end{proof}
In \cite{ATY} we proved Theorem \ref{type2asymp} for $n=2$ using a different approach from the present one.

From this theorem the characterization of type $1$ representations follows just as in the one-variable case.  We obtain a strengthening of Theorem \ref{thm1.2}.
\begin{theorem}\label{type1asymp}
 The following three conditions are equivalent for a function $h\in\pick$.
\begin{enumerate}
\item  The function $h$ has a Nevanlinna representation of type $1$;
\item 
\[
\liminf_{s\to\infty} s\abs{h(is\chi)}<\infty;
\]
\item 
\beq 
\label{type1lim} \lim_{s\to\infty} s\abs{h(is\chi)}<\infty.
\eeq
\end{enumerate}
\end{theorem}
\begin{proof}
  We follow Lax's treatment \cite{lax02} of the one-variable Nevanlinna theorem.

(1)$\Rightarrow$(3)  Suppose that $h$ has a type $1$ representation as in equation \eqref{type1formula} for some $\h,\ A,\ Y$ and $v$.  Then 
\begin{align}
h(is\chi) &= \ip{(A - is)^{-1}\alpha}{\alpha} \notag \\
	    &= \ip{(A + is)(A^2 + s^2)\inv \alpha}{\alpha} \notag,
\end{align}
and so
\[
\re sh(is\chi) = \ip{sA (A^2 + s^2)\inv \alpha}{\alpha}, \quad \im sh(is\chi) = \ip{s^2 (A^2 + s^2)\inv \alpha}{\alpha}.
\]
Let $\nu$ be the scalar spectral measure for $A$ corresponding to the vector $\al\in\h$.  Then
\[
\re sh(is\chi) = \int \frac{st}{t^2+s^2}\ \dd\nu(t), \quad \im sh(is\chi) = \int \frac{s^2}{t^2+s^2} \ \dd\nu(t).
\]
The integrand in the first integral tends pointwise in $t$ to $0$ as $s\to\infty$, and by the inequality of the means it  is no greater than $\half$; thus the integral tends to $0$ as $s\to\infty$ by the Dominated Convergence Theorem.  The integrand in the second integral increases monotonically to $1$ as $s\to \infty$.  Thus
\[
\re sh(is\chi) \to 0, \qquad \im sh(is\chi) \to \|\al\|^2 \quad \mbox{ as } s\to \infty.
\]
Hence the inequality \eqref{type1lim} holds.  Thus (1)$\Rightarrow$(3).

(3)$\Rightarrow$(2) is trivial.

(2)$\Rightarrow$(1) Suppose that 
\beq\label{toomany}
\liminf_{s \to \infty} s\abs{h(is\chi)} < \infty. 
\eeq
As
\beq \notag
\liminf_{s \to \infty} s \IM h(is\chi) \leq \liminf_{s \to \infty} s\abs{h(is\chi)} < \infty,
\eeq
$h$ satisfies condition \eqref{type2liminf} of Theorem \ref{type2asymp}. Therefore $h$ has a representation of type $2$, say
\[
 h(z) = a + \ip{(A - z_Y)^{-1}\alpha}{\alpha}.
\]
It remains to show that $a = 0$. The inequality \eqref{toomany} implies that there exists a sequence $s_n$ tending to $\infty$ such that $h(is_n\chi) \to 0$. But
\beq \notag
 \RE h(is_n\chi) = a + \ip{A(A^2 + s_n^2)\inv\alpha}{\alpha} \to a.
\eeq
Hence $a = 0$ and $h$ has a type $1$ representation. This establishes (2)$\Rightarrow$(1).
\end{proof}

\section{Carapoints at infinity}\label{caraInfty}

How can we recognise from function-theoretic properties whether a given function in the $n$-variable Loewner class admits a Nevanlinna representation of a given type?  In the preceding section it was shown that it depends on growth along a single ray through the origin. In this section we describe the notion of carapoints at infinity for a function in the Pick class, and in the next section we shall give succinct criteria for the four types in the language of carapoints.  

Carapoints (though not with this nomenclature) were first introduced by Carath\'eodory  in  1929 \cite{car29} for a function $\ph$ on the unit disc, as a hypothesis in the ``Julia-Carath\'eodory Lemma".   For any $\tau\in\T$, a function $\vp$ in the Schur class {\em satisfies the Carath\'eodory condition at $\tau$} if
\beq\label{julia}
\liminf_{\lambda\to\tau} \frac{1 - \abs{\vp(\lambda)}}{1 - \abs{\lambda}} < \infty.
\eeq
The notion has been generalized to other domains by many authors.  Consider domains $U\subset\C^n$ and $V\subset \C^m$ and an analytic function $\ph$ from $U$ to the closure of $V$.   The function $\ph$  is said to satisfy Carath\'eodory's condition at $\tau\in\partial U$ if 
\[
\liminf _{\la\to\tau}\frac{ \dist(\ph(\la), \partial V)}{\dist(\la, \partial U)} < \infty.
\]
Thus, for example,  when $U=\Pi^n, V=\Pi$,  a function $h\in \Pick_n$ satisfies Carath\'eodory's condition at the point $x\in\R^n$ if
\beq \label{caragen}
\liminf_{z\to x} \frac {\im h(z)}{\min_j \im z_j } < \infty.
\eeq
This definition works well for finite points in $\partial U$, but for our present purpose we need to consider points at infinity in the boundaries of $\Pi^n$ and $\Pi$.   We shall introduce a variant of Carath\'eodory's condition for the class $\Pick_n$ with the aid of the Cayley transform
 \beq\label{cayleypar}
  z = i\frac{1+\lambda}{1-\lambda}, \qquad \la = \frac{z-i}{z+i},
 \eeq
which furnishes a conformal map between $\D$ and $\Pi$, and hence a biholomorphic map between $\D^n$ and $\Pi^n$ by co-ordinatewise action.  We obtain a one-to-one correspondence between $\schur\setminus \{\bf {1}\}$ and $\Pick_n$ via the formulae
 \beq \label{cayleyfunc}
  h(z) = i\frac{1+\vp(\lambda)}{1-\vp(\lambda)},  \quad \ph(\la)= \frac{h(z)-i}{h(z)+i}
 \eeq
where $\bf{1}$ is the constant function equal to $1$ and $\la, z$ are related by equations \eqref{cayleypar}.  For $\ph\in\schur$ we define $\tau\in\T^n$ to be a {\em carapoint} of $\ph$ if
\beq\label{juliaBis}
\liminf_{\lambda\to\tau} \frac{1 - |\vp(\lambda)|}{1 - \norm{\lambda}_\infty} < \infty.
\eeq
We can now extend the notion of carapoints to points at infinity.  The point $(\infty,\dots,\infty)$ in the boundary of $\Pi^n$ corresponds to the point $\chi$ in the closed unit disc;
 as in the last section, $\chi$ denotes the point $(1,\dots,1)\in\C^n$.  

 \begin{definition}
  Let $h$ be a  function in the Pick class $\Pick_n$ with associated function $\vp$ in the Schur class $\schur$ given by equation \eqref{cayleyfunc}.  Let $\tau \in \T^n, \ x \in (\R\cup\infty)^n$ be related by
\beq\label{xandtau}
x_j=i\frac{1+\tau_j}{1-\tau_j} \quad \mbox{ for } j=1,\dots,n.
\eeq
 We say that $x$ is a {\em carapoint} for $h$ if $\tau$ is a carapoint for $\vp$.   We say that $h$ has a {\em carapoint at} $\infty$ if $h$ has a carapoint at $(\infty,\dots,\infty)$, that is, if $\vp$ has a carapoint at $\chi $.
 \end{definition}
Note that, for a point $x\in\R^n$, to say that $x$ is a carapoint of $h$ is {\em not} the same as saying that $h$ satisfies the Carath\'eodory condition \eqref{caragen} at $x$.  Consider the function $h(z)=-1/z_1$ in $\Pick_n$.  Clearly $h$ does not satisfy Carath\'eodory's condition at $0\in\R^n$. However, the function $\ph$ in $\schur$ corresponding to $h$ is $\ph(\la) = -\la_1$, which does have a carapoint at $-\chi$, the point in $\T^n$ corresponding to $0 \in\R^n$.  Hence $h$ has a carapoint at $0$.  

We shall be mainly concerned with carapoints at $0$ and $\infty$.  The following observation will help us identify them.  For any $h\in\Pick_n$ we define $h^\flat\in\Pick_n$ by
\[
h^\flat(z)= h\left(-\frac{1}{z_1},\dots, -\frac{1}{z_n}\right)\quad\mbox{ for }z\in\Pi^n.
\]
For $\ph\in\schur$ we define 
\[
\ph^\flat(\la)= \ph(-\la).
\]
If $h$ and $\ph$ are corresponding functions,  as in equations \eqref{cayleyfunc}, then so are $h^\flat$ and $\ph^\flat$.
\begin{proposition}\label{caracriter}
The following conditions are equivalent for a function $h\in\Pick_n$.
\begin{enumerate}
\item  $\infty$ is a carapoint for $h$;
\item $0$ is a carapoint for $h^\flat$;
\item
\[
\liminf_{y\to 0+} \frac{\im h^\flat(iy\chi)}{y|h^\flat(iy\chi)+i|^2} < \infty;
\]
\item
\[
\liminf_{y\to\infty} \frac{y\im h(iy\chi)}{|h(iy\chi)+i|^2} < \infty.
\]
\end{enumerate}
\end{proposition}
\begin{proof}
(1)$\Leftrightarrow$(2)
Since $-\chi\in\T^n$ corresponds under the Cayley transform to $0\in \R^n$, we have
\begin{align*}
\infty \mbox{ is a carapoint of  } h\quad  &\Leftrightarrow \quad  \chi \mbox{  is a carapoint of } \ph  \\
	&\Leftrightarrow \quad  -\chi \mbox{ is a carapoint of }\ph^\flat  \\
	&\Leftrightarrow   \quad 0  \mbox{ is a carapoint of } h^\flat.
\end{align*}

\nin (2)$\Leftrightarrow$(3)
A consequence of the $n$-variable Julia-Carath\'eodory Theorem \cite{jafari,abate}, 
is that $\tau\in\T^n$ is a carapoint of $\ph\in\schur$ if and only if
\[
\liminf_{r\to 1-} \frac{1-|\ph(r\tau)|}{1-r} < \infty.
\]
It follows that
\begin{align*}
0 \mbox{ is a carapoint for } h^\flat \quad &\Leftrightarrow \quad  -\chi \mbox{ is a carapoint for } \ph^\flat \\
	&\Leftrightarrow \quad \liminf_{r\to 1-}  \frac{1-|\ph^\flat(-r\chi)|}{1-r} < \infty \\
	&\Leftrightarrow \quad \liminf_{r\to 1-}  \frac{1-|\ph^\flat(-r,-r)|^2}{1-r^2} < \infty. \\
\end{align*}
Let $iy\in\Pi$ be the Cayley transform of $-r \in (-1,0)$, so that $y\to 0+$ as $r\to 1-$.
In view of  the identity
\beq\label{caracond}
\frac{1-|\ph(\la)|^2}{1-\|\la\|^2_\infty}= \left(\max_j \frac{|z_j+i|^2}{\im z_j}\right) \frac{\im h(z)}{|h(z)+i|^2}
\eeq
we have
\begin{align*}
0 \mbox{ is a carapoint for } h^\flat \quad &\Leftrightarrow \quad \liminf_{y\to 0+} \frac{|iy+i|^2}{y} \frac{\im h^\flat(iy\chi)}{|h^\flat(iy\chi)+i|^2} < \infty \\
	&\Leftrightarrow \quad \liminf_{y\to 0+} \frac{\im h^\flat(iy\chi)}{y|h^\flat(iy\chi)+i|^2} < \infty.
\end{align*}
\nin (3)$\Leftrightarrow$(4)  Replace $y$ by $1/y$.
\end{proof}

\begin{corollary}\label{suffcond}
If $f\in\Pick_n$ satisfies Carath\'eodory's condition
\beq\label{ccond}
\liminf_{z\to x} \frac{\im f(z)}{\im z} < \infty
\eeq
at $x\in\R^n$ then $x$ is a carapoint for $f$.  If
\[
\liminf_{y\to\infty} y\im f(iy\chi) < \infty
\]
then $\infty$ is a carapoint for $f$.
\end{corollary}
\begin{proof}
Let $h=f^\flat \in\Pick_n$.  Clearly $|h^\flat(z)+i| \geq 1$ for all $z\in\Pi^n$.  If the condition \eqref{ccond} holds for $x=0$ then
\[
\liminf_{z\to 0} \frac{\im  h^\flat(z)}{|h^\flat(z)+i|^2 \min_j \im z_j}\leq \liminf_{z\to 0}\frac {\im h^\flat(z)}{\min_j\im z_j}< \infty
\]
and hence, by (2)$\Leftrightarrow$(3) of Proposition \ref{caracriter}, $0$ is a carapoint for $h^\flat=f$.  The case of a general $x\in\R^n$ follows by translation.
\end{proof}
If $h \in \Pick_n$ has a carapoint at $x\in(\R\cup\infty)^n$ then it has a value at $x$ in a natural sense. If $\vp\in\schur$ has a carapoint at $\tau\in\T^n$, then by \cite{jafari}  there exists a unimodular constant $\vp(\tau)$ such that 
 \beq\label{ntlimit}
  \lim_{\lambda \nt \tau} \vp(\lambda) = \vp(\tau).
 \eeq 
 Here $\la\nt\tau$ means that $\la$ tends nontangentially to $\tau$ in $\D^n$. 

 \begin{definition}
  If $h\in\Pick_n$ has a carapoint at $x\in (\R\cup\infty)^n$ then we define
  \[
   h(x) = \threepartdef{\infty}{\vp(\tau)=1}{}{}{\ds i\frac{1+\vp(\tau)}{1-\vp(\tau)}}{\vp(\tau)\neq 1}
  \]
where $\tau\in\T^n$ corresponds to $x$ as in equation \eqref{xandtau}.
 \end{definition}
Thus $h(\infty)\in \R\cup\{\infty\}$ when $\infty$ is a carapoint of $h$.

In the example $h(z)=-1/z_1$, since the value of $\ph(-\la)$ at $-\chi$ is $1$, we have $h(0)=\infty$.

Although the value of $h(\infty)$ is defined in terms of the Schur class function $\ph$, it can be expressed more directly in terms of  $h$. 

\begin{proposition} \label{hinfty}
 If $\infty$ is a carapoint of $h$ then
\beq\label{ntlimatinfty}
h(\infty)= h^\flat(0) = \lim_{z \nt \infty} h(z).
\eeq
\end{proposition}
Here we say that $z \nt \infty$ if $z\to (\infty,...,\infty)$ in the set $\{z\in\Pi^n: (-1/z_1,\dots, -1/z_n) \in S\}$ for some set $S\subset \Pi^n$ that approaches $0$ nontangentially, or equivalently, if $z\to (\infty,\dots,\infty)$ in a set on which $\|z\|_\infty /\min_j \im z_j$ is bounded.
\begin{proof}
Clearly
\[
h(\infty) =\infty  \quad \Leftrightarrow  \quad \ph(\chi)=1  \quad \Leftrightarrow \quad \ph^\flat(-\chi)=1 \quad \Leftrightarrow  \quad h^\flat(0) =\infty.
\]
Similarly, for $\xi\in\R$,
\[
h(\infty)=\xi  \quad \Leftrightarrow \quad\ph(\chi) = \frac{\xi -i}{\xi+i} \quad\Leftrightarrow  \quad \ph^\flat(-\chi) = \frac{\xi -i}{\xi+i} \quad \Leftrightarrow \quad  h^\flat(0) =\xi.
\]
Thus, whether $h(\infty)$ is finite or infinite, $h(\infty)=h^\flat(0)$.
Equation \eqref{ntlimatinfty} follows from the relation \eqref{ntlimit}.
\end{proof}

\section{Types of functions in the Loewner class}\label{caraTypes}
In this section we shall show that the type of a function $h\in\p$ is entirely determined by whether or not $\infty$ is a carapoint of $h$ and by the value of $h(\infty)$.
Let us make precise the notion of the {\em type} of a function in $\p$.
\begin{definition}
  A function $h \in\pick$ is of {\em type $1$} if it has a Nevanlinna representation of type $1$.   For $n=2, 3$ or $4$ we say that $h$ is of {\em type $n$} if $h$ has a Nevanlinna representation of type $n$ but has no representation of type $n-1$.
 \end{definition}
Clearly every function in $\pick$ is of exactly one of the types $1$ to $4$.  We shall now prove Theorem \ref{typeoffcn}.  Recall that it states the following, for any function $h \in \p$.
\begin{enumerate}
\item $h$ is of type $1$ if and only if $\infty$ is a carapoint of $h$ and $h(\infty) = 0$;
 \item $h$ is of type $2$ if and only if $\infty$ is a carapoint of $h$ and $h(\infty)\in \R\setminus\{0\}$;
 \item $h$ is of type $3$ if and only if $\infty$ is not a carapoint of $h$;
 \item $h$ is of type $4$ if and only if $\infty$ is a carapoint of $h$ and $h(\infty) = \infty$.
\end{enumerate}
\begin{proof}
(2)  Let $h\in\p$ have a type $2$ representation $h(z)=a+\ip{(A-z_Y)\inv v}{v}$ with $a\neq 0$.  By Theorem \ref{type2asymp},
\[
\liminf_{y\to\infty} y \im h(iy\chi) < \infty.
\]
By Corollary \ref{suffcond}, $\infty$ is a carapoint for $h$.  Furthermore, by Proposition \ref{hinfty}
\[
h(\infty) =\lim_{y\to\infty} h(iy\chi) = a \in\R\setminus\{0\}.
\]

Conversely, suppose that $\infty$ is a carapoint for $h$ and $h(\infty)\in \R\setminus\{0\}$.  By Proposition \ref{caracriter}
\[
\liminf_{y\to\infty}\frac{y\im h(iy\chi)}{|h(iy\chi)+i|^2} < \infty
\]
while by Proposition \ref{hinfty}
\[
\lim_{y\to\infty} |h(iy\chi)+i|^2 = h(\infty)^2+1 \in (1,\infty).
\]
On combining these two limits we find that
\[
\liminf_{y\to\infty} y\im h(iy\chi) < \infty,
\]
and so, by Theorem \ref{type2asymp}, $h$ has a representation of type $2$.  Since $h(\infty)\neq 0$ it is clear that $h$ does not have a representation of type $1$.  Thus (2) holds.

A trivial modification of the above argument proves that (1) is also true.\\

\nin (4)  Let $h$ be of type 4. Then $h$ has no type $3$ representation, and so, by Theorem \ref{type3asymp}, there exists $\delta > 0$ and a sequence $(s_n)$ of positive numbers tending to $ \infty$ such that 
\[
 \frac{1}{s_n} \IM h(is_n\chi) \geq \delta > 0.
\]
Let $y_n=1/s_n$; then $-1/(is_n) = iy_n$, and we have
\beq\label{hflatyn}
 y_n \im h^\flat(iy_n\chi) \geq \delta \quad \mbox{ for all } n\geq 1.
\eeq
Since $|h^\flat(z)+i| > \im h^\flat (z)$ for all $z$, we have
\begin{align*}
\liminf_{z\to 0} \frac{\im h^\flat(z)}{|h^\flat(z)+i|^2 \min_j \im z_j} & \leq \liminf_{z\to 0} \frac{1}{ \im h^\flat(z) \min_j \im z_j} \\
	&\leq \liminf_{n\to\infty} \frac{1}{y_n \im h^\flat(iy_n\chi)} \\
	& \leq 1/\delta.
\end{align*}
Hence $(0,0)$ is a carapoint of $h^\flat$, and so $\infty$ is a carapoint of $h$.

Since $y_n \to 0$ it follows from the inequality \eqref{hflatyn} that $\im h^\flat(iy_n\chi) \to \infty$, hence that $h^\flat(0) =\infty$, and therefore that $h(\infty)=\infty$.

Conversely, suppose that $\infty$ is a carapoint of $h$ and that $h(\infty) = \infty$. We shall show that
 \beq \label{nolimit}
\lim_{s\to\infty}\frac{1}{s}\IM h(is\chi) \neq 0,
\eeq
and it will follow from Theorem \ref{type3asymp} that $h$ does not have a representation of type $3$, that is, $h$ is of type $4$.

Let $\ph\in\schur$ correspond to $h$ and let $r\in (0,1)$ correspond to $is\in\Pi$.  Then
\begin{align}\label{need1}
\frac{1}{s}\im h(is\chi) &= \frac{1-r}{1+r} \quad  \frac{1-|\ph(r\chi)|^2}{|1-\ph(r\chi)|^2}  \nn \\
	&= \frac{1-|\ph(r\chi)|^2}{1-r^2} \, \, \frac{(1-r)^2}{|1-\ph(r\chi)|^2}.
\end{align}
By hypothesis, $\chi$ is a carapoint for $\ph$ and $\ph(\chi)=1$.  By definition of carapoint,
\[
\liminf_{z\to\chi} \frac{1-|\ph(z)|^2}{1-\|z\|_\infty^2} = \al < \infty \quad \mbox{ for all } s > 0.
\]
The $n$-variable Julia-Carath\'eodory Lemma (see \cite{jafari,abate}) 
 now tells us that $\al > 0$ and 
\beq\label{need2}
\frac{|1-\ph(r\chi)|^2}{|1-r|^2} \leq  \al \frac{1-|\ph(r\chi)|^2}{1-r^2} \quad \mbox{ for all } r\in (0,1).
\eeq
On combining equations \eqref{need1} and \eqref{need2} we obtain
\[
\frac{1}{s}\im h(is\chi) \geq \frac{1}{\al} > 0 \quad \mbox{ for all } s>0.
\]
Thus the relation \eqref{nolimit} is true, and so, by Theorem \ref{type3asymp},  $h$ is of type $4$.

Statement (3) now follows easily.  The function $h\in\p$ is of type $3$ if and only if it is not of types $1,2$ or $4$, hence if and only if it is not the case that $\infty$ is a carapoint for $h$ and $h(\infty) \in \R\cup\{\infty\}$, hence if and only if $\infty$ is not a carapoint of $h$.
\end{proof}

We now show that there are functions in the Pick class $\Pick_2$ of all four types. We return to Example \ref{4types} and show that the functions in $\Pick_2$ which we presented there are indeed of the stated types.
\begin{example}\label{fourtypes} \rm
(1) The function
\[
h(z)= -\frac{1}{z_1+z_2} = \ip{(0-z_Y)\inv v}{v}_{\C},
\]
where $Y=\half$ and $v=1/\sqrt{2}$,
is obviously of type $1$.  Let us nevertheless check that $\infty$ is a carapoint of $h$ and $h(\infty) = 0$, in accordance with Theorem \ref{typeoffcn}.  We have $h(iy,iy) = \half i/y$
and hence 
\[
\liminf_{y\to 0+} y\im h(iy,iy)= \half.
\]
  Thus $\infty$ is a carapoint for $h$ by Proposition \ref{caracriter}.   Moreover $h(iy,iy) \to 0$ as $y\to \infty$, and therefore $h(\infty)=0$.\\

\nin (2) It is immediate that the function $1+h$, with $h$ as in (1), is of type 2, and that $\infty$ is a carapoint of $1+h$ with value $1$.\\

\nin (3)  We have seen that the function
\beq\label{truehbis}
h(z) =\threepartdef{ \ds  \frac{1}{1+z_1z_2}\left(z_1-z_2 + \frac{iz_2(1+z_1^2)}{\sqrt{z_1z_2}}\right)}{z_1z_2 \neq -1}{}{}{\half (z_1+z_2)}{z_1z_2 = -1}
\eeq
has a representation of type $3$. To show that $h$ is indeed of type $3$ we must prove that $\infty$ is not a carapoint of $h$.

For all $y>0$ we have $h(iy,iy)=i$.  Hence
\[
\liminf_{y\to \infty} \frac{y\im h(iy,iy)}{|h(iy,iy)+i|^2} = \liminf_{y\to\infty} \frac{y}{4} = \infty.
\]
By Proposition \ref{caracriter}, $\infty$ is not a carapoint for $h$.  Thus $h$ is of type $3$.\\

\nin (4)  The function
\[
h(z)= \frac{z_1z_2}{z_1+z_2} = -1 \left/ \left(-\frac{1}{z_1} - \frac{1}{z_2}\right) \right.
\]
is clearly in $\Pick_2$.  We gave a type $4$ representation of $h$ in Example \ref{4types}.
We claim that $\infty$ is a carapoint of $h$.  We have $h(iy,iy)= \half iy$, and thus
\begin{align*}
\liminf_{y\to \infty} \frac{y \im h(iy,iy)}{|h(iy,iy)+i|^2} &= \liminf_{y\to \infty} \frac{\half y^2}{|\half iy+i|^2}=2.
\end{align*}
Hence $\infty$ is a carapoint for $h$.  Furthermore $h(iy,iy)= \half iy \to \infty$ as $y\to \infty$, and so $h(\infty)=\infty$.  Thus $h$ is of type $4$.
\end{example}

Another example of a function of type $4$ is $h(z)=\sqrt{z_1z_2}$.

\section{ Rates of growth in the Loewner class } \label{growth}
The Nevanlinna representation formulae give rise to growth estimates for functions in the $n$-variable Loewner class.    It turns out that growth is mild, both at infinity and close to the real axis.  Even though the type of a function is determined by its growth on the single ray $\{iy\chi:y>0\}$, in turn the growth of the function on the entire polyhalfplane is constrained by its type.

Consider first the one-variable case.  If $h$ is the Cauchy transform of a finite positive measure $\mu$ then 
\[
|h(z)| \leq \int \frac{\dd \mu(t)}{|t-z|} \leq \int \frac{\dd\mu(t)}{\im z}  = \frac{C}{\im z}
\]
for some $C>0$ and for all $z\in\Pi$.  For a general function $h$ in the Pick class, by Nevanlinna's representation (Theorem \ref {thm1.1}) there exist $a\in\R, b\geq 0$ and a finite positive measure $\mu$ on $\R$ such that, for all $z\in\Pi$,
\begin{align*}
h(z) &= a + bz + \int \frac{1+tz}{t-z} \ \dd\mu(t)\\
	&= a+bz +\int \frac{1+z^2}{t-z} + z \ \dd\mu(t)
\end{align*}
and therefore
\begin{align*}
|h(z)| &\leq |a| + b|z| + \left(\frac{1+|z|^2}{\im z} + |z|\right) \mu(\R)\\
	& \leq C\left(1+ |z| + \frac{1+|z|^2}{\im z}\right)
\end{align*}
for some $C>0$.

Similar estimates hold for the Loewner class.
\begin{proposition}
For any function $h\in\p$ there exists a non-negative number $C$ such that, for all $z\in\Pi^n$,
\beq\label{type4gr}
|h(z)| \leq C\left( 1+\|z\|_1 + \frac{1+\|z\|_1^2}{\min_j \im z_j}\right).
\eeq

For any function $h\in\p$ of type $2$ there exists a non-negative number $C$ such that, for all $z\in\Pi^n$,
\beq\label{type2gr}
|h(z)| \leq C\left(1+ \frac{1}{\min_j\im z_j}\right).
\eeq

For any function $h\in\p$ of type $1$ there exists a non-negative number $C$ such that, for all $z\in\Pi^n$,
\beq\label{type1gr}
|h(z)| \leq  \frac{C}{\min_j\im z_j}.
\eeq
\end{proposition}
\begin{proof}
Let $h\in\p$.  Let $\N,\M, A, P, a$ and $v$ be as in Theorem \ref{thm2.4}, so that
\[
h(z)= a+\ip{M(z)v}{v}
\]
for all $z\in\Pi^n$, where $M(z)$ is the matricial resolvent given by equation \eqref{defM}.  By Proposition \ref{2x2resolv} we have, for all $z\in\Pi^n$,
\begin{align*}
\|M(z)\| & \leq (1+\sqrt{10} \|z\|_1)\left(1 + \frac{1+\sqrt{2}\|z\|_1}{\min_j \im z_j}\right) \\
	&\leq 1+\sqrt{10} \|z\|_1 + B \frac{1+\|z\|_1+\|z\|_1^2}{\min_j \im z_j}
\end{align*}
for a suitable choice of $B \geq 0$.  Hence 
\begin{align*}
|h(z)| &\leq |a| + \|M(z)\| \|v\|^2 \\
	&\leq |a| + \left( 1+\sqrt{10}\|z\|_1 + B \frac{1+\|z\|_1+\|z\|_1^2}{\min_j \im z_j}\right) \|v\|^2.
\end{align*}
Since
\[
1+\|z\|_1+\|z\|_1^2 \leq \tfrac 32 (1+\|z_1\|^2),
\]
we have
\[
|h(z)| \leq C\left( 1+\|z\|_1 + \frac{1+\|z\|_1^2}{\min_j \im z_j}\right)
\]
for some choice of $C>0$ and for all $z\in\Pi^n$.  Thus the estimate \eqref{type4gr} holds.

Similarly, the estimates \eqref{type2gr} and \eqref{type1gr} follow easily from the simple resolvent estimate \eqref{boundAzYinv}.
\end{proof}

\section{ Structured resolvent identities }\label{resolvident}
To conclude the paper we point out that there are structured analogs of the classical resolvent identity
\[
(A-z)\inv -(A-w)\inv = (z-w)(A-z)\inv (A-w)\inv
\]
 for any $z,w$ in the resolvent set of an operator $A$.

\begin{proposition} \label{resIdent}
Let $A$ be a densely defined self-adjoint operator on a Hilbert space $\h$ and let $Y$ be a positive decomposition of $\h$.  For all $z,w \in \Pi^n$
\beq\label{resId2}
(A-z_Y)\inv - (A-w_Y)\inv= (A-z_Y)\inv(z-w)_Y(A-w_Y)\inv.
\eeq
If $M(z)$ is the structured resolvent of type $3$ corresponding to $A$ and $Y$ then
\beq\label{resId3}
M(z)-M(w)\left| \mathcal D(A) \right. =(1-iA)(A-z_Y)\inv(z-w)_Y(A-w_Y)\inv(1+iA).
\eeq
\end{proposition}
\begin{proof}
The first of these identities is immediate.  For the second, by equation \eqref{3rdtype3},
\begin{align*}
M(z)-M(w)\left| \mathcal D(A) \right. =(1-iA)\left((A-z_Y)\inv-(A-w_Y)\inv\right)(1+iA),
\end{align*}
and the  identity \eqref{resId3} follows from \eqref{resId2}.
\end{proof}
\begin{proposition}\label{10.2}
Let $\mathcal H$ be the orthogonal direct sum of Hilbert spaces $\mathcal{ N, M}$, let $A$ be a densely defined self-adjoint operator on $\mathcal M$ with domain $\mathcal{D}(A) $ and let $P$ be an orthogonal decomposition of $\mathcal H$.  For every $z, \ w\in\Pi^n$, as operators on $\N\oplus\mathcal{D}(A)$,
\begin{align}\label{resId4}
M(z)-M(w) &= \bbm -i&0\\0&1-iA \ebm\left(\bbm1&0\\0&A\ebm -z_P\bbm 0&1\\0&0\ebm\right)\inv (z-w)_P\nn \\
	&\hspace{2cm} \times \left(\bbm 1&0\\0&A\ebm - \bbm 0&0\\0&1\ebm w_P\right)\inv\bbm i&0\\0&1+iA \ebm.
\end{align}
\end{proposition}
\begin{proof}
Let 
\[
D= \bbm i&0\\0&1+iA \ebm : \N\oplus \mathcal{D}(A) \to \h.
\]
By equations \eqref{simpler} and \eqref{simpler2} we have
\begin{align*}
M(z)&-M(w)\left| \N\oplus\mathcal{D}(A)\right.\\
	& = D^*\left\{\left(\bbm1&0\\0&A\ebm -z_P\bbm 0&1\\0&0\ebm\right)\inv\left(z_P\bbm 0&0\\0&1\ebm+ \bbm 0&0 \\ 0&1\ebm \right)\right. \\
	&\hspace*{1cm}\left. -\left(\bbm 1&0\\0&0\ebm w_P + \bbm 0&0\\0&1\ebm\right)\left(\bbm 1&0\\0&A\ebm - \bbm 0&0\\0&1\ebm w_P\right)\inv \right\} D \\
	& = D^*\left(\bbm1&0\\0&A\ebm -z_P\bbm 0&1\\0&0\ebm\right)\inv \left\{\left(z_P\bbm 0&0\\0&1\ebm+ \bbm 0&0 \\ 0&1\ebm \right)\left(\bbm 1&0\\0&A\ebm  - \bbm 0&0\\0&1\ebm w_P\right)\right. \\
	&\left. \hspace*{0.3cm} - \left(\bbm1&0\\0&A\ebm -z_P\bbm 0&1\\0&0\ebm\right)\left(\bbm 1&0\\0&0\ebm w_P + \bbm 0&0\\0&1\ebm\right)  \right\} \left(\bbm 1&0\\0&A\ebm - \bbm 0&0\\0&1\ebm w_P\right)\inv D.
\end{align*}
The term in braces in the last expression reduces to $(z-w)_P$, and the identity \eqref{resId4} follows.
\end{proof}
\begin{corollary}
With the assumptions of Proposition {\rm \ref{10.2}}, there exists an analytic function $F:\Pi^n\to \mathcal{L}(\h)$ such that, for all $z,\ w\in\Pi^n$,
\beq\label{MFP}
M(z)- M(w) = F(\bar z)^*(z-w)_P F(w).
\eeq
\end{corollary}
The statement follows from Proposition \ref{10.2} just as Proposition \ref{Mloewner} follows from Proposition \ref{alternM}.  If $F$ is defined by equation \eqref{defFz} then $F(z)$ is a bounded operator on $\h$, $F$ is analytic on $\Pi^n$ and Proposition \ref{10.2} states that equation \eqref{MFP} holds  on $\N\oplus \mathcal{D}(A)$.  It follows by continuity that equation \eqref{MFP} holds  on $\h$.

\bibliographystyle{plain}

\nin J. Agler, Department of Mathematics, University of California at San Diego, CA 92103, USA.\\

\nin R. Tully-Doyle,  Department of Mathematics, University of California at San Diego, CA 92103, USA.\\

\nin N. J. Young, School of Mathematics, Leeds University, Leeds LS2 9JT {\em and} School of Mathematics and Statistics, Newcastle University, Newcastle upon Tyne NE3 4LR, England.  Email N.J.Young@leeds.ac.uk

\end{document}